\documentclass[a4paper,reqno]{amsart}
\usepackage{amsmath,amssymb}
\usepackage{enumerate}
\usepackage{xspace}
\usepackage{boxedminipage}
\usepackage{listings}
\usepackage{tabularx}
\usepackage[font=footnotesize]{subfig}
\usepackage{algpseudocode}
\usepackage{tristan}
\usepackage{fullpage}

\numberwithin{equation}{section}
\setboolean{showtodo}{true}
\setboolean{shownotes}{true}
\setboolean{showchanges}{false}
\setboolean{usemathrsfs}{true}
\author{
  Elizabeth Mansfield
}
\address{
  Elizabeth Mansfield \newline
  School of Mathematics, Statistics \& Actuarial Science\newline
  University of Kent\newline
  Canterbury\newline
  GB-CT2 7NF , England UK}
\email{\linkedurl{E.L.Mansfield@kent.ac.uk}}

\author{
  Tristan Pryer
} 
\address{
  Tristan Pryer \newline
  Department of Mathematics and Statistics\newline
  Whiteknights\newline
  PO Box 220\newline
  Reading\newline
  GB-RG6 6AX , England UK}
\email{\linkedurl{T.Pryer@reading.ac.uk}}
\thanks{The authors were supported by the EPSRC grant EP/H024018/1}

\title[Noether type FE conserved quantities]{{Noether type
  discrete conserved quantities arising from a Finite Element
  approximation of a variational problem}}

\date{\today}
\pdfformat{true}
\renewcommand{\highlight}[1]{#1}
\begin{document}
\keywords{finite element method, conserved quantities, Noether's
  theorem, variational problem.}

\subjclass{65N303; 49M25; 22E99}

\maketitle
\begin{abstract}
  In this work we prove a weak Noether type theorem for a class of
  variational problems which include \emph{broken extremals}. We then
  use this result to prove discrete Noether type conservation laws for
  certain classes of finite element discretisation of a model elliptic
  problem. In addition we study how well the finite element scheme
  satisfies the continuous conservation laws arising from the
  application of Noether's 1st Theorem (E. Noether 1918).

  We summarise extensive numerical tests, illustrating the conservativity
  of the discrete Noether law using the $p$--Laplacian as an example.
\end{abstract}

\section{Introduction \highlight{and historical background}}
\label{sec:intro}

\highlight
{The purpose of this paper is to show that variational numerical
problems have their own conservation laws which derive from the same
principle as that discovered by Noether, giving rise to discrete
(numerical) forms of conservation laws which are automatically
preserved by the scheme.}

Symmetries are an extremely important and continually occurring
feature of differential equations arising from many applicable areas,
including mathematical physics, meteorology and differential geometry,
that was first developed by Sophus Lie for the purpose of studying
solutions of differential equations in the late 19th century
\cite[c.f.]{Lie:1971}.

Noether's (1st) Theorem \cite{Noether:1971} is a striking result which
in the continuous setting connects these symmetries with conservation
laws associated to the Euler--Lagrange equations {of a variational
problem}. Roughly, the theorem states that given a variational problem
with an underlying symmetry, there exists a natural conservation law
associated to it. For example, a symmetry of translation with respect
to the spatial coordinates results in a conservation of linear
momentum, a symmetry of rotation results in conservation of angular
momentum and a symmetry of translation with respect to the temporal
coordinate gives a conservation of energy. A famous example from
meteorology is that of \emph{potential vorticity}. This is a
conservation law arising from a particle relabelling pseudo-group
symmetry. This quantity is extremely important in studying the
evolution of a cyclone \highlight{\cite[c.f.]{DavisEmanuel:1991}}.

The work of Noether has gained public attention recently with the
publication of an article in the New York Times \cite{Angier:2012}
where the result is
\begin{quote}
  ``consider[ed] \dots as important as Einstein's theory of
  relativity; it underg\highlight{[ir]}ds much of today's vanguard research in
  physics, including the hunt for the almighty Higgs boson.''
\end{quote}

In the discrete setting, Noether's Theorem has been studied in terms
of difference equations \cite{Dorodnitsyn:2001,HydonMansfield:2004},
where it was shown that a discrete equivalent of the conservation law
holds when a discrete symmetry was built into the discrete
Lagrangian. In this work we turn our attention to the finite element
method (FEM). FEMs form one of the most successful numerical methods
for approximating the solution to partial differential equations
(PDEs) \cite[c.f.]{Braess:2001, Brenner:1994, Ciarlet:1978}. A topic
which has been the subject of much ongoing research is that of
constructing FEMs which inherit some property of the continuous
problem. The notion of discretisations inheriting some geometric
property from the continuous problem can be seen as a generalisation
of geometric integration \cite[c.f.]{HairerLubichWanner:2006} to the
case of PDEs and is a rapidly developing area of research. Some of the
properties studied in the discretisation of PDEs are the same as in
the geometric integration of the ODE, for example the Hamiltonian
structure of a given problem. Others are based on completely new
notions, for example, the recent development of the discrete exterior
calculus \cite{ArnoldFalkWinther:2010,RobidouxSteinberg:2011}, which,
as the name suggests, is a discrete equivalent to the Cartan based
exterior calculus. This has allowed for a rigorous description of
discrete differential forms and the associated discrete function
spaces as a discrete differential complex. This provides a framework
which may be used as a first step in the construction of a variational
complex in a similar light to that developed in
\cite{HydonMansfield:2004} for difference equations. A first step in
this direction was taken in \cite{Mansfield:2006}. A review of some of
the huge quantity of topics arising from this area, including Lie
group integrators, discrete gradient methods as well as FEMs for
differential forms is given in
\cite{ChristiansenMunthe-KaasOwren:2011}.

As opposed to geometric integrators, the term used for numerical
methods with some geometric property of an ODE, the methods for PDEs
are generally called \emph{mimetic methods}. The class of FEMs which
fall under the mimetic framework are the mixed methods, for example
the Raviart--Thomas scheme \cite{RaviartThomas:1977}. It is not only
FEMs which fall under the category of mimetic methods, in fact there
are finite difference (FD) \cite{BrezziBuffaLipnikov:2009} and finite
volume (FV) schemes which are characterised as mimetic. Note that also
there is an intrinsic relationship between each of them. For example,
an appropriate choice of quadrature for the Raviart-Thomas finite
element scheme results in the mimetic FD scheme
\cite{CangianiManziniRusso:2009}.

The Lagrangian piecewise polynomial FEM is not a mimetic method. Most
standard methods cannot inherit geometric properties of the continuous
PDE. There is an underlying algebraic condition which must be
satisfied for these properties to be inherited by the approximation
scheme \cite{HydonMansfield:2004,Mansfield:2006}.


The classical Noether theorem is only applicable to classical
solutions of the variational problem. As such we derive weaker
versions of the theorem applicable to a wider class of solutions to
the problem, including the \emph{broken extremals}. We will discuss
how these laws are naturally passed down to the Lagrangian finite
element scheme and hence quantify the discrete Noether quantity
associated to this FEM. That is, we write the exact Noether quantity
for this discretisation, in the same spirit as
\cite{HydonMansfield:2004}.

We will also study how well the Lagrangian finite element scheme
satisfies the strong conserved quantities arising from Noether's
Theorem measured in an appropriate weak norm. That is, we consider how
well this finite element scheme approximates the Noether conservation
law for the continuous problem (when one exists). We will also present
some interesting numerical results, quantifying the deviation of the
approximation in terms of a computable estimator which we are able to
use to construct an $h$--adaptive scheme (local mesh refinement) aimed
at minimising the violation of the smooth conservation law to a user
specified tolerance.

The paper is set out as follows: In \S\ref{sec:notation} we introduce
some fundamental notation and the model problem we consider. In
\S\ref{sec:noether} we briefly describe Noether's Theorem and the
background material needed. To illustrate its application we apply the
Theorem to a simple model problem. In \S\ref{sec:noether-weak} we weaken the
invariance criterion on which the classical Noether Theorem is based,
ultimately allowing us prove two versions of the theorem applicable to
weaker solutions of the problem. In \S\ref{sec:discretisation} we
discuss how the results of \S\ref{sec:noether-weak} can be passed down
to give discrete counterparts to our weak Noether's Theorem. We perform
numerical experiments to demonstrate that the quantities derived are
indeed conserved at the discrete level. We also discuss trivial Lie
group actions (those of translation with respect to the dependent
variable) and how the mimetic methods relate to this case. Finally, in
\S\ref{sec:numerics-lag} we study the properties of the finite element
solution with respect to the original (strong) Noether Theorem. We
also detail an interesting numerical result by constructing a
computable estimator, aimed at measuring the violation of the strong
Noether Theorem in a weak norm for the Lagrangian finite element
scheme. We perform some numerical experiments demonstrating that there
is a \emph{superconvergence} of the estimator over the finite element
approximation of the solution to the Euler--Lagrange equations. We
then proceed to test an adaptive scheme based on the estimate allowing
us to minimise the discrete violation of the continuous conserved
quantity up to user specified tolerance.

\section{Notation}
\label{sec:notation}

Let $\W \subset\reals^d$ be a bounded domain with
boundary $\partial \W$. We begin by introducing the Sobolev spaces
\cite{Ciarlet:1978,Evans:1998}
\begin{gather}
  \leb{p}(\W)
  =
  \ensemble{\phi}
           {\int_\W \norm{\phi}^p \highlight{\d \geovec x} < \infty}   \text{ for } p\in[1,\infty) 
             \text{ and }
             \leb{\infty}(\W)
             =
             \ensemble{\phi}
                      {\esssup_\W \norm{\phi} < \infty},
                      \\
                      \sob{k}{p}(\W) 
                      = 
                      \ensemble{\phi\in\leb{p}(\W)}
                               {\partial^{\vec\alpha}\phi\in\leb{p}(\W), \text{ for } \norm{\geovec\alpha}\leq k}
                               \text{ and }
    \sobh{k}(\W)
    := 
    \sob{k}{2}(\W),
\end{gather}
which are equipped with the following norms and semi-norms:
\begin{gather}
  \Norm{v}_{\leb{p}(\W)}^p
  :=
       {\int_\W \norm{v}^p \highlight{\d \geovec x}},
       \qquad
       \Norm{v}_{\sob{k}{p}(\W)}^p
       = 
       \sum_{\norm{\vec \alpha}\leq k}\Norm{\partial^{\vec \alpha} v}_{\leb{p}(\W)}^p,
       \\
       \norm{v}_{\sob{k}{p}(\W)}^p
       =
       \sum_{\norm{\vec \alpha} = k}\Norm{\partial^{\vec \alpha} v}_{\leb{p}(\W)}^p,
       \qquad 
       \Norm{v}_{\sobh{k}(\W)}^2 = \Norm{v}_{\sob{k}{2}(\W)}^2,
\end{gather}
where $\vec\alpha = \{ \alpha_1,\dots,\alpha_d\}$ is a
multi-index, $\norm{\vec\alpha} = \sum_{i=1}^d\alpha_i$ and
derivatives $\partial^{\vec\alpha}$ are understood in a weak sense. We pay
particular attention to the cases $k = 1,2$ and
\begin{gather}
  \sobz{1}{p}(\W) := \text{closure of }\cont{\infty}_0(\W) \text{ in } \sob{1}{p}(\W).
\end{gather}
Let $L = L\qp{\geovec x, u, \nabla u}$ be the \emph{Lagrangian}. We
will let
\begin{equation}
  \label{eq:action-functional}
  \dfunkmapsto[.]
  {\cJ[\cdot]}
  \phi
  {\sobz{1}{p}(\W)}
  {\cJ[\phi] := \displaystyle\int_\W L(\geovec x, \phi, \nabla \phi) \highlight{\d \geovec x}}
  {\reals}
\end{equation}
be known as the \emph{action functional}. The problem arising from the
calculus of variations is to seek a function \emph{extremising} the
action functional. For simplicity we will consider the minimisation
problem, {that is, to find} $u\in\sobz{1}{p}(\W)$ such that
\begin{equation}
  \label{eq:variational-prob}
  \cJ[u] = \inf_{v\in\sobz{1}{p}(\W)} \cJ[v].
\end{equation}
{Note that we are implicitly coupling the minimisation problem with
homogeneous Dirichlet boundary conditions.}

We will use the notation that 
\highlight{\begin{equation}
  \partial_1 q := \nabla q
  =
  \Transpose{\qp{\frac{\partial q\qp{\geovec x, u, \nabla u}}{\partial x_1}
      ,
      \dots
      ,
      \frac{\partial q\qp{\geovec x, u, \nabla u}}{\partial x_d} }
  }
\end{equation}}
denotes a column vector of spatial derivatives of \highlight{a generic scalar valued function $q$}, \ie
derivatives with respect to the independent variables. The derivative
with respect to the dependent variable is denoted 
\highlight{\begin{equation}
\partial_2 q := \frac{\partial
  q\qp{\geovec x, u, \nabla u}}{\partial u}
\end{equation}}
 and let $\geovec p = \Transpose{\qp{p_1,\dots,p_d}}
= \nabla u$ then
\highlight{
\begin{equation}
  \partial_3 q := 
  \Transpose{\qp{\frac{\partial q\qp{\geovec x, u, \nabla u}}{\partial p_1}
      ,
      \dots
      ,
      \frac{\partial q\qp{\geovec x, u, \nabla u}}{\partial p_d}}
  }
\end{equation}}
denotes the vector of derivatives of $q$ with respect to the gradient
of $u$ componentwise. \highlight{We use $\div$ to represent the
  \emph{spatial divergence} of a vector valued function, $\geovec q = \qp{q_1, \dots, q_d}$, hence
  \begin{equation}
    \partial_1 \geovec q := \div\qp{\geovec q}
    =
    \frac{\partial q_1\qp{\geovec x, u, \nabla u}}{\partial x_1}
    +
    \dots
    +
    \frac{\partial q_d\qp{\geovec x, u, \nabla u}}{\partial x_d}.
  \end{equation}
  The derivative with respect to the independent variable is then a
  column vector
  \begin{equation}
    \partial_2 \geovec q 
    :=
    \Transpose{\qp{
        \frac{\partial q_1\qp{\geovec x, u, \nabla u}}{\partial u}
        ,\dots ,
        \frac{\partial q_d\qp{\geovec x, u, \nabla u}}{\partial u}
      }}
  \end{equation}
  and
  \begin{equation}
    \partial_3 \geovec q
    :=
    \begin{bmatrix}
      \frac{\partial q_1\qp{\geovec x, u, \nabla u}}{\partial p_1}
      ,\dots ,
      \frac{\partial q_d\qp{\geovec x, u, \nabla u}}{\partial p_1}
      \\
      \vdots\qquad \ddots \qquad\vdots
      \\
      \frac{\partial q_1\qp{\geovec x, u, \nabla u}}{\partial p_d}
      ,\dots ,
      \frac{\partial q_d\qp{\geovec x, u, \nabla u}}{\partial p_d}
      \end{bmatrix}.
  \end{equation}
  With the above notations we may
  introduce the \emph{total derivative operator}, defined for scalar
  valued functions as
  \begin{equation}
    \label{eq:total-deriv}
    \D q\qp{\geovec x, u, \nabla u} 
    :=
    \partial_1 q
    +
    \nabla u \partial_2 q
    +
    \Hess u \partial_3 q,
  \end{equation}
  and \emph{total divergence operator}, for vector valued functions as
  \begin{equation}
    \label{eq:total-div}
    \Div\qp{\geovec q\qp{\geovec x, u,\nabla u}}
    :=
    \partial_1 \geovec q
    +
    \Transpose{\qp{\partial_2 \geovec q}} \nabla u 
    +
    \frob{\partial_3 \geovec q}{\Hess u},
  \end{equation}
  where $\frob{\geomat X}{\geomat Y} = \trace{\Transpose{\geomat
      X}{\geomat Y}}$ denotes the \emph{Frobenious inner product}
  between matrices.  }

It is well known \cite[c.f.]{Evans:1998, GiaquintaHildebrandt1:1996}
that if $u$ is a (smooth) minimiser of the variational problem
(\ref{eq:variational-prob}) then it solves the quasilinear, second
order PDE called the \emph{Euler--Lagrange equations}
\begin{equation}
  \label{eq:Euler-Lagrange-equations}
  \EL[u] :=
  -
  \Div\qp{\partial_3 L} 
  +
  \partial_2 L = 0.
\end{equation}


\section{Noether's First Theorem}
\label{sec:noether}

For the reader's benefit we will briefly describe Noether's first
theorem in the continuous, smooth case and necessary background
material. We assume, in this section, that $L$ is smooth and the
minimisation problem \eqref{eq:variational-prob} has a solution (not
necessarily unique) {which is at least $\cont{2}(\W)$, \ie
  smooth enough to satisfy the Euler--Lagrange equations \eqref{eq:Euler-Lagrange-equations}}

\begin{Defn}[one-parameter group]
  \label{defn:1-parameter-group}
  The transformation 
  \begin{equation}
    \label{eq:transformation}
    \qp{\geovec x, u} \to \qp{\Xi(\geovec x, u; \epsilon), \Phi(\geovec x, u; \epsilon)} 
    =:
    \qp{\widetilde{\geovec x}, \widetilde u}
  \end{equation}
  is said to be a \emph{one-parameter group} if the following conditions hold
  \begin{enumerate}

  \item The parameter choice of $\epsilon = 0$ yields the identity,
    \ie 
    \begin{equation*}
      \qp{\geovec x,u} = \qp{\Xi(\geovec x, u; 0), \Phi(\geovec x, u; 0)}.
    \end{equation*}

  \item The inverse is given by the parameter $-\epsilon$, \ie
    \begin{equation*}
      \qp{\geovec x,u} 
      =
      \qp{\Xi(\widetilde{\geovec x}, \widetilde u; -\epsilon), \Phi(\widetilde{\geovec x}, \widetilde u; -\epsilon)}.
    \end{equation*}

  \item The transformation is closed under composition, \ie if
    \begin{equation*}
      \qp{\widehat{\geovec x}, \widehat u}
      =
      \qp{\Xi(\widetilde{\geovec x}, \widetilde u; \delta), \Phi(\widetilde{\geovec x}, \widetilde u; \delta)}
    \end{equation*}
    then 
    \begin{equation*}
      \qp{\widehat{\geovec x}, \widehat u}
      =
      \qp{\Xi({\geovec x}, u; \epsilon+\delta), \Phi({\geovec x}, u; \epsilon+\delta)}.
    \end{equation*}
  \end{enumerate}
\end{Defn}

\begin{Defn}[infinitesimal]
  \label{defn:infinitesimal}
  The infinitesimals, $\geovec \xi(\geovec x,u)$ and $\phi(\geovec
  x,u)$ of the one parameter group are defined as
  \begin{gather}
    \geovec \xi(\geovec x,u) 
    := 
    \lim_{\epsilon\to 0}
    \frac{\d \Xi(\geovec x, u; \epsilon)}{\d \epsilon}
    \\
    \phi(\geovec x,u) 
    := 
    \lim_{\epsilon \to 0}\frac{\d \Phi(\geovec x, u; \epsilon)}{\d \epsilon}
  \end{gather}
\end{Defn}

\begin{Defn}[characteristics]
\label{defn:characteristics}
  {We define the characteristics, which are given in terms of the
  infinitesimals of the group, to be}
  \begin{equation}
    Q\qp{\geovec x, u, \nabla u}
    := 
    \phi(\geovec x, u) - \Transpose{\qp{\geovec \xi(\geovec x, u)}}\nabla u.
  \end{equation}
\end{Defn}

\begin{Defn}[variational symmetry]
  \label{defn:var-symm}
  Let $\Gamma := \setof{\qp{\geovec x, u(\geovec x)} : \geovec x \in
    \Upsilon}$ be the graph of $u$ over a subdomain such that
  $\closure\Upsilon \subset \W$. Also let $\Upsilon_{\Xi} =
  \Xi\qp{\Gamma; \epsilon}$, then the transformation
  (\ref{eq:transformation}) is said to be a variational
  symmetry if
{
  \begin{equation}
    \int_{\Upsilon} L\qp{\geovec x, u, \nabla u} \d \geovec x
    =
    \int_{\Upsilon_\Xi} L\qp{\widetilde{\geovec x}, \widetilde u, \widetilde{\nabla u}} \d \widetilde {\geovec x}
  \end{equation}
}
  holds for any smooth subdomain $\Upsilon$ of $\W$.
\end{Defn}

\begin{The}[infinitesimal invariance {\cite[Thm 4.12]{Olver:1993}}]
  \label{the:inf-inv}
  A variational symmetry group with infinitesimals $\geovec \xi, \phi$
  and characteristics $Q$ of the action functional
  \begin{equation}
    \cJ[u] = \int_\W L(\geovec x, u, \nabla u) \highlight{\d \geovec x}
  \end{equation}
  satisfies
  \begin{equation}
    0 
    =
    \Transpose{\qp{\D Q}}
    \partial_3 L
    +
    Q 
    \partial_2 L
    + 
    \Div\qp{L \geovec\xi}
  \end{equation}
\end{The}
\begin{Proof}
  See {\cite[Thm 4.12]{Olver:1993}} 
\end{Proof}

\begin{The}[Noether's 1st Theorem {\cite[Thm 4.29]{Olver:1993}}]
  \label{the:noether}
  Suppose the variational problem (\ref{eq:variational-prob}) is
  invariant under the action of a one-parameter group of symmetries
  with characteristics $Q$. Then $Q$ is also a characteristic of a
  conservation law of the Euler--Lagrange equation
  \eqref{eq:Euler-Lagrange-equations}, that is, there exists a vector valued
  functional $\cC = \cC[u]$ such that
  {
  \begin{equation}
    \label{eq:conservation-laws}
    \Div\qp{\cC[u]} = Q \EL[u].
  \end{equation}
  Hence over solutions of the Euler--Lagrange equations $\EL[u] =
  0$, we have that
  \begin{equation}
    \Div\qp{\cC[u]} = 0.
  \end{equation}
  For the problem we consider in this work, that of a first order
  Lagrangian, the conservation law, $\cC$, takes the form}
  \begin{equation}
    \label{eq:noether-functional}
    \cC[u] 
    = 
    -
    \qp{
      L \geovec \xi 
      +
      {\partial_3 L \qp{\phi - \Transpose{\geovec \xi}\nabla u}}}.
  \end{equation}
\end{The}
\begin{Proof}
  Using the result of Theorem \ref{the:inf-inv} we have that
  \begin{equation}
    \begin{split}
      0 
      &=
      \Transpose{\qp{\D Q}}
      \partial_3 L
      +
      Q 
      \partial_2 L
      + 
      \Div\qp{L \geovec\xi}.
    \end{split}
  \end{equation}
  Noting by the product rule that
  \begin{equation}
    \Transpose{\qp{\D Q}}
    \partial_3 L
    =
    -Q \Div{\partial_3 L}
    +
    \Div\qp{Q \partial_3 L}
  \end{equation}
  then it holds that
  \begin{equation}
    \begin{split}  
    \label{eq:EL-cons-law}
      0 
      &=
      - 
      Q \Div\qp{\partial_3 L}
      +
      \Div\qp{Q \partial_3 L}
      +
      Q 
      \partial_2 L
      + 
      \highlight{\Div\qp{L \geovec\xi}}
      \\
      &=
      Q \qp{-\Div\qp{\partial_3 L} + \partial_2 L}
      + 
      \Div\qp{Q \partial_3 L + L \geovec\xi}.
    \end{split}
  \end{equation}
  This concludes the proof with
  \begin{equation}
    \cC[u] = \highlight{-}\qp{Q \partial_3 L + L \geovec\xi},
  \end{equation}
  as required. 
\end{Proof}

\begin{Rem}[the form of $\cC$]
  \label{rem:form-of-c}
  It is clear from the identity (\ref{eq:conservation-laws}) that for
  our model problem, \highlight{that of minimising a first order variational
  problem \eqref{eq:variational-prob}, we have} $\cC = \cC(u, \nabla
  u)$.
\end{Rem}

\begin{Rem}[the beauty of the theorem]
  What makes Theorem \ref{the:noether} truly remarkable is its
  constructive nature. For completeness we will give an example of the
  construction of $\cC$ for the Laplacian.
\end{Rem}

\begin{Example}[Laplace's problem]
  \label{ex:laplace}
  Let us consider the case $f = f(\norm{\geovec x})$ then the
  Lagrangian,
  \begin{equation}
    L\qp{\geovec x, u, \nabla u} := \frac{1}{2} \norm{\nabla u}^2 - f u,
  \end{equation}
  is invariant under the rotational group $SO(d)$. For simplicity we
  restrict \highlight{to the case} $d=2$, set $\geovec x =
  \Transpose{\qp{x, y}}$, then we calculate the infinitesimals from
  the group of rotations, note that in this case $\Phi \equiv 0$ and
  \begin{equation}
    \geovec \Xi(\geovec x, u;\epsilon) 
    = 
    \begin{bmatrix}
      x\cos\epsilon - y\sin\epsilon
      \\
      x\sin\epsilon \highlight{+} y\cos\epsilon
    \end{bmatrix}.
  \end{equation}
  It then holds that
  \begin{equation}
    \lim_{\epsilon\to 0}\frac{\d \geovec \Xi(\geovec x;\epsilon)}{\d \epsilon}
    =
    \begin{bmatrix}
      -y
      \\
      x
    \end{bmatrix}.
  \end{equation}
  In this case the characteristic of the group of rotations is
  \begin{equation}
    y \partial_x u - x \partial_y u.
  \end{equation}
  Making use of Theorem \ref{the:noether} we see
  \begin{equation}
    \cC[u] 
    =
    \begin{bmatrix}
      y \qp{\qp{\partial_y u}^2 - \qp{\partial_x u}^2}/2
      + x \partial_x u \partial_y u + y f u
      \\
      x \qp{\qp{\partial_y u}^2 - \qp{\partial_x u}^2}/2
      - y \partial_x u \partial_y u - x f u
    \end{bmatrix}
  \end{equation}
  is a conservation law over solutions of $\EL[u] = 0$.
\end{Example}

\begin{Rem}[trivial Lie group actions]
  \label{rem:trivial-actions}
  For any variational problem, the Euler--Lagrange equations, as
  already mentioned, are given in variational (or divergence) form. As
  such, if we assume that $L$ does not depend on $u$, that is $L =
  L(\geovec x, \nabla u)$, then the Euler--Lagrange equations
  themselves \highlight{are} a Noether conservation law. Indeed, consider the
  case of Example \ref{ex:laplace} with $f\equiv 0$. It is clear by
  definition that $\Delta u = \div\qp{\nabla u} = 0$ is a conservation
  law. It arises from Noether's Theorem under the trivial Lie group
  action, that of translation in the dependent variable
  \begin{equation}
    \qp{\geovec x, u} \to \qp{\geovec x, u + \epsilon}.
  \end{equation}
  For this action, the infinitesimals are $\geovec \xi = \geovec 0$ and
  $\phi = 1$.
\end{Rem}

\section{Noether's theorem for weak solutions}
\label{sec:noether-weak}

Noether's Theorem (Theorem \ref{the:noether}) as it is stated in
\S\ref{sec:noether} only makes sense for classical solutions of the
Euler--Lagrange equations (\ref{eq:Euler-Lagrange-equations}). We wish
to ``weaken'' the theorem such that it is applicable to {extremals}
which are $\sob{1}{\infty}(\W)$ with jump discontinuities \highlight{in the derivatives}, the so
called \emph{broken extremals} \cite{GiaquintaHildebrandt1:1996}. 
We begin by defining a weaker invariance condition than that of
Definition \ref{defn:var-symm}.
\begin{Defn}[weak variational symmetry]
  \label{def:weak-var-symm}
  The transformation (\ref{eq:transformation}) is said to be a weak
  variational symmetry if 
  {
  \begin{equation}
    \label{eq:weak-var-symm}
    \int_{\W} L\qp{\geovec x, u, \nabla u} \d \geovec x
    =
    \int_{\W_\Xi} L\qp{\widetilde{\geovec x}, \widetilde u, \widetilde{\nabla  u}} \d \widetilde {\geovec x}
  \end{equation}}
  holds over the domain $\W$.  
\end{Defn}

\begin{Rem}[strong symmetry $\implies$ weak symmetry]
  We note that any strong variational symmetry is also a weak symmetry
  but the converse is not true.
\end{Rem}

\begin{The}[Noether type conserved quantities for weak variational symmetries]
  \label{the:noether-weak}
  Suppose that the variational problem \eqref{eq:variational-prob} has
  a weak variational symmetry.
  Let $\phi$ and $\geovec \xi$ be the infinitesimal generators of the
  symmetry as in Definition \ref{defn:infinitesimal}. Then
  \begin{equation}
    \begin{split}
      \label{eq:weak-cons-law}
      0 &= \int_\W 
      \Transpose{ \qp{ \partial_3 L }} 
      \D \phi
      +
      \partial_2 L \phi 
      +
      L \Div\qp{\geovec \xi} 
      +
      \Transpose{\qp{\partial_1 L}}\geovec \xi
      -
      \Transpose{\qp{ \partial_3 L}}\nabla
      u \Div\qp{\geovec \xi} 
      \highlight{\d \geovec x}.
    \end{split}
  \end{equation}
\highlight{  Over smooth minimisers, \ie if $u\in\cont{2}(\W)$. We have
  \begin{equation}
    \begin{split}
      \label{eq:weak-cons-law2}
      0 
      &=
      \int_\W \Div\qp{L \geovec \xi} - \Div\qp{\partial_3 L} Q +
      \partial_2 L Q 
      +
      \int_{\partial \W} \Transpose{ \qp{Q \partial_3
          L}}\geovec n \d s
      \\
      &=
      \int_{\partial \W} \Transpose{ \qp{ Q\partial_3
          L + \xi L}}\geovec n \d s
    \end{split}
  \end{equation}}
\end{The}

\begin{Rem}[structure of (\ref{eq:weak-cons-law})]
  The weak conservation law given in Theorem \ref{the:noether-weak}
  has a very clear structure. The first two terms \highlight{in
    (\ref{eq:weak-cons-law})} represent the weak Euler--Lagrange
  equations
  . The last three terms in
  (\ref{eq:weak-cons-law}) represent the weak conservation
  law itself.
\end{Rem}

{ 
\begin{Proof}[\highlight{of Theorem \ref{the:noether-weak}}]
  Using the fact that the problem (\ref{eq:variational-prob}) has a
  weak variational symmetry, from Definition \ref{def:weak-var-symm}
  we see that
  \begin{equation}
    0
    =
    \lim_{\epsilon \to 0}
    \frac{1}{\epsilon}
    \qp{
      \int_{\W_\Xi} L\qp{\widetilde{\geovec x}, \widetilde u, \widetilde{\nabla  u}} \d \widetilde {\geovec x}
      -
      \int_{\W} L\qp{\geovec x, u, \nabla u} \d \geovec x 
    }.
    \end{equation}
    Using a coordinate transformation from $\W_\Xi$ to $\W$, \highlight{we have}
    \begin{equation}
      0 
      =
      \lim_{\epsilon \to 0}
      \frac{1}{\epsilon}
      \qp{
        \int_{\W} L\qp{\widetilde{\geovec x}, \widetilde u, \widetilde{\nabla  u}} \frac{\d \widetilde {\geovec x}}{\d \geovec x} \d \geovec x
        -
        \int_{\W} L\qp{\geovec x, u, \nabla u} \d \geovec x 
      },
    \end{equation}
    \highlight{noting that $\widetilde{x}=\widetilde{x}(x,u(x))$, so that the first integrand is indeed defined on $\Omega$.}  Making use of Definition \ref{defn:infinitesimal} \highlight{ and the fact that 
    \begin{equation}
      \widetilde{\nabla u} = \nabla u + \epsilon \qp{\D \phi - \nabla u \Div{\geovec \xi}} + \Oh(\epsilon^2),
    \end{equation}}
  it holds that
  \begin{equation}
    \label{eq:inf-var-sym}
    0 
    =
    \int_\W
    \Transpose{\qp{\partial_1 L}}\geovec \xi
    + 
    \partial_2 L \phi
    +
    \Transpose{\qp{\partial_3 L}} \D \phi
    -
    \Transpose{\qp{\partial_3 L}} \nabla u \Div\qp{\geovec \xi}
    +
    L\Div\qp{\geovec \xi}
    \d \geovec x,
  \end{equation}
  as required for the first equality. The second arises from noting
  (\ref{eq:inf-var-sym}) implies
  \begin{equation}
    \begin{split}
      0 
    &=
    \int_\W
    \Transpose{\qp{\D L - \nabla u\partial_2 L - \Hess u\partial_3 L}}\geovec \xi
    + 
    \partial_2 L \phi
    +
    \Transpose{\qp{\partial_3 L}} \D \phi
    \\
    &\qquad -
    \Transpose{\qp{\partial_3 L}} \nabla u \Div\qp{\geovec \xi}
    +
    L\Div\qp{\geovec \xi}
    \d \geovec x
    \\
    &=
    \int_\W 
    \Div\qp{L\geovec \xi} 
    +
    \pd2 L \qp{\phi - \Transpose{\qp{\nabla u}}\geovec \xi}
    +
    \Transpose{\qp{\pd3 L}} \qp{\D \phi - \nabla u \Div{\geovec \xi} - \Hess u \geovec \xi}
    \d \geovec x.
    \end{split}
  \end{equation}
  Using the fact that
  \begin{equation}
    \D Q 
    =
    \D \phi - \nabla u \Div\qp{\geovec \xi} - \Hess u \geovec \xi
  \end{equation}
  we have that
  \begin{equation}
    \begin{split}
      0 
      &=
      \int_\W \Div\qp{L\geovec \xi}
      +
      Q \partial_2 L
      - Q \Div\qp{\partial_3 L}
      +
      \int_{\partial{\W}} Q \Transpose{\qp{\partial_3 L}}\geovec n
    \end{split}
  \end{equation}
  upon applying Stokes Theorem and noting that $u$ is now an extremal
  hence satisfies the Euler--Lagrange equations (modulo natural
  boundary conditions).

  Note that if the group action is separable we may separate the proof
  into computing the \emph{inner} and \emph{outer} variations with
  respect to the infinitesimals of the one parameter group (see
  Definition \ref{defn:infinitesimal})
  \cite[c.f.]{GiaquintaHildebrandt1:1996}, where the inner variations
  are with respect to the independent variables and the outer
  variation with respect to the dependent variables.

  \end{Proof}
}

\begin{Cor}[strong conservation law $\implies$ weak conservation law]
  Let the variational problem \eqref{eq:variational-prob} have a
  variational symmetry in the sense of Definition \ref{defn:var-symm}
  and that the minimiser to the variational problem is smooth
  $u\in\cont{2}(\W)$, then \eqref{eq:weak-cons-law} holds.
\end{Cor}

Now we have developed the framework sufficiently to state our main
result in this section. Here we are concerned with \emph{broken
  extremals}, that is, functions whose derivatives have finitely many
jump discontinuities.

\begin{Defn}[broken extremal]
  An extremal, $u\in\cont{0}(\W)$, to the problem
  (\ref{eq:variational-prob}) is said to be a broken extremal if it is
  piecewise $\cont{2}(\W)$ over the domain $\W$ with bounded
  derivative. {That is, $\W$ can be decomposed into finitely many
  open subsets, $\{\W_i\}_{i=1}^N$ such that 
  \begin{enumerate}
  \item the subsets make up the entire domain, \ie $\W = \union{i}\closure{\W_i}$,
  \item they are non-overlapping, \ie $\W_i\cap \W_j = \emptyset$ and
  \item the solution is smooth over each of the subsets, \ie
    $u\in\cont{2}(\W_i)\cap \sob{1}{\infty}(\W)$.
  \end{enumerate}}
\end{Defn}

\begin{Defn}[skeleton and jumps]
  We define
  \begin{equation}
    \cF := \union{i} \setof{\geovec x \in \partial \W_i}
  \end{equation}
  to be the \emph{skeleton} of the decomposition. \highlight{We will
    assume that the domain is decomposed in such a way that the
    skeleton is Lipschitz continuous.} Let $\geovec n_i$ be the
  outward pointing normal to $\W_i$, we then define \emph{jumps} of
  scalars and vector valued functions as
  \begin{gather}
    \jump{v} := v|_{\W_1} \geovec n_1 + v|_{\W_2} \geovec n_2
    \\
    \jump{\geovec v} 
    := 
    \Transpose{\qp{\geovec v|_{\W_1}}}  \geovec n_1 
    +
    \Transpose{\qp{\geovec v|_{\W_2}}} \geovec n_2,
  \end{gather}
  respectively. 
\end{Defn}

\begin{Defn}[piecewise variational symmetry]
  \label{def:pw-var-sym}
  {Let $u$ be a broken extremal to the variational problem
    (\ref{eq:variational-prob}) and let $\setof{\W_i}_{i=1}^N$ be the
    decomposed domain of $u$. Then \eqref{eq:transformation} is a
    piecewise variational symmetry if it is a variational symmetry
    over $\W_i$ for each $i=1,\dots N$.}
\end{Defn}

\begin{The}[conserved quantities for broken extremals of the variational problem]
  \label{the:noether-broken}
  Suppose the variational problem (\ref{eq:variational-prob}) has a
  piecewise variational symmetry. Then 
  \begin{equation}
    \begin{split}
      0 &=
      \sum_i \int_{\W_i} \qp{-\Div\qp{\pd 3 L} + \pd2 L }\phi - L
      \pd1{\geovec \xi} - \Transpose{\qp{\partial_3 L}} \nabla u
      \Div\qp{\geovec \xi} + L\Div\qp{\geovec \xi} \d \geovec x 
      \\
      &      \qquad\qquad\qquad
      \qquad\qquad\qquad
 +
      \int_{\cF} \jump{L \Transpose{\geovec \xi} +
        \Transpose{\qp{\pd3 L}} \phi }\d s.
    \end{split}
  \end{equation}
  Over broken extremals we have that
  \begin{equation}
    \begin{split}
      0
      &=
      \sum_{i=1}^N\int_{\W_i}
      L \Div{\geovec \xi}
      -
      L \pd1 {\geovec \xi}
      -
      \Transpose{\qp{\pd 3 L}}
      \nabla u
      \Div{\geovec \xi} \highlight{\d \geovec x}
      \\
      &
      \qquad\qquad\qquad
      \qquad\qquad\qquad
      +
      \int_{\cF} \jump{\pd 3 L \phi - L\geovec \xi} \d s.
    \end{split}
  \end{equation}
\end{The}
\highlight
{\begin{Proof}
  The proof of this result follows the same lines as the proof of
  Theorem \ref{the:noether-weak}. Using \eqref{eq:inf-var-sym} we have
  that for each $\W_i$
  \begin{equation}
    0 
    =
    \duality{\Transpose{\qp{\partial_1 L}}}{\geovec \xi}
    +
    \int_{\W_i}
    \partial_2 L \phi
    +
    \Transpose{\qp{\partial_3 L}} \D \phi
    -
    \Transpose{\qp{\partial_3 L}} \nabla u \Div\qp{\geovec \xi}
    +
    L\Div\qp{\geovec \xi}
    \d \geovec x,
  \end{equation}
  where we use $\duality{\cdot}{\cdot}$ to denote the duality action
  on a Sobolev space from its dual. Hence we have that the first term
  is understood in a duality sense, \ie
  \begin{equation}
    \duality{\Transpose{\qp{\pd 1 L}}}{\geovec \xi}_{\W_i}
    =
    -\int_{\W_i}
    L \div\qp{\geovec \xi} \d \geovec x
    + \int_{\partial\W_i} L\Transpose{\geovec \xi } \geovec n \d s.
  \end{equation}
  Integrating by parts we see
  \begin{equation}
    \begin{split}
      0 &= \int_{\W_i} \qp{-\Div\qp{\pd 3 L} + \pd2 L }\phi - L
      \pd1{\geovec \xi} - \Transpose{\qp{\partial_3 L}} \nabla u
      \Div\qp{\geovec \xi} + L\Div\qp{\geovec \xi} \d \geovec x 
      \\
      &\qquad +
      \int_{\pd{} \W_i} L \Transpose{\geovec \xi} \geovec n +
      \Transpose{\qp{\pd3 L}} \geovec n \phi \d s.
    \end{split}
  \end{equation}
  Summing over each of the subdomains we have
  \begin{equation}
    \begin{split}
      0 &=
      \sum_i \int_{\W_i} \qp{-\Div\qp{\pd 3 L} + \pd2 L }\phi - L
      \pd1{\geovec \xi} - \Transpose{\qp{\partial_3 L}} \nabla u
      \Div\qp{\geovec \xi} + L\Div\qp{\geovec \xi} \d \geovec x 
      \\
      &\qquad +
      \int_{\cF} \jump{L \Transpose{\geovec \xi} +
      \Transpose{\qp{\pd3 L}} \phi }\d s,
    \end{split}
  \end{equation}
  as required for the first equality. For the second we note that over
  broken extremals the Euler--Lagrange equations vanish over each
  $\W_i$, concluding the proof.
\end{Proof}
}

\section{Finite element conservation laws}
\label{sec:discretisation}

\subsection{Discretisation}
In this section we calculate the discrete counterpart to Theorem
\ref{the:noether-broken} in the finite element context. To that end,
let $\T{}$ be a conforming triangulation of $\W$, namely, $\T{}$ is a
finite family of sets such that
\begin{enumerate}
\item $K\in\T{}$ implies $K$ is an open simplex (segment for $d=1$,
  triangle for $d=2$, tetrahedron for $d=3$),
\item for any $K,J\in\T{}$ we have that $\closure K\meet\closure J$ is
  a full subsimplex (i.e., it is either $\emptyset$, a vertex, an
  edge, a face, or the whole of $\closure K$ and $\closure J$) of both
  $\closure K$ and $\closure J$ and
\item $\union{K\in\T{}}\closure K=\closure\W$.
\end{enumerate}

We let $\E{}$ be the skeleton (set of \highlight{internal} common
interfaces) of the triangulation $\T{}$ and say $e\in\E$ if $e$ is on
the interior of $\W$ and $e\in\partial\W$ if $e$ lies on the boundary
$\partial\W$.

The \emph{shape regularity} of $\T{}$ is defined as
\begin{equation}
  \label{eqn:def:shape-regularity}
  \mu(\T{}) := \inf_{K\in\T{}} \frac{\rho_K}{h_K},
\end{equation}
where $\rho_K$ is the radius of the largest ball contained inside
$K$ and $h_K$ is the diameter of $K$.
 We use the convention where $\funk h\W\reals$ denotes the
 \emph{meshsize function} of $\T{}$, i.e.,
 \begin{equation}
   h(\geovec{x}):=\max_{\closure K\ni \vec x}h_K,
 \end{equation}
 where $h_K$ is the diameter of an element K. We introduce the
 \emph{finite element spaces}
 \begin{gather}
   \label{eqn:def:finite-element-space}
   \fes
   :=
   \ensemble{\Phi \in \cont{0}(\W)
   }{\Phi\vert_{K} \in \poly k\Forall K\in\T{}}
   \\
   \feszero = \fes \cap \hoz(\W),
 \end{gather}
 where $\poly k$ denotes the linear space of polynomials in $d$
 variables of degree no higher than a positive integer $k$. We
 consider $k\geq 1$ to be fixed and denote by $N := \dim{\feszero}$.

The Galerkin approximation to the variational problem
(\ref{eq:variational-prob}) is to seek $U\in\feszero\subset\hoz(\W)$ such that
\begin{equation}
  \label{eq:dis-min-prob}
  \cJ[U] = \inf_{V\in\feszero} \cJ[V]
\end{equation}

The finite element scheme defined by (\ref{eq:dis-min-prob}) is
guaranteed to be well posed under some assumptions on $L$ that allow
us to invoke the Lax--Milgram Theorem or the more generally applicable
inf-sup condition \cite{ErnGuermond:2004}. Henceforth we now assume the
continuous minimisation problem admits a unique solution.

We may now proceed in deriving a finite element Noether type conservation
law. As already seen, the conservation law arises after taking inner
and outer variations of the variational problem. The outer variation
can be characterised by the following Lemma.

\begin{Lem}[discrete Euler--Lagrange equations]
  \label{lem:dis-EL-eqns}
  The discrete Euler--Lagrange equations associated to the variational
  minimisation problem \eqref{eq:dis-min-prob} are to seek $U\in\feszero$
  such that
  \begin{equation}
    \begin{split}
      0
      &=
      \int_\W 
      \qp{-\Div_K\qp{\pd3 L}
        +
        \pd 2 L } V \highlight{\d \geovec x}
      +
      \int_\E
      \jump{\pd3 L} V 
      \highlight{\d \geovec s}
      \Foreach V\in\feszero,
    \end{split}
  \end{equation}
  where $L=L(\geovec x, U, \nabla U)$.
\end{Lem}
\begin{Proof}
  Define the real valued function which we call the \emph{outer
    variation operator}
  \begin{equation}
    o\qp{\epsilon}
    :=
    \int_\W 
    L\qp{\geovec x, U + \epsilon V, \nabla U + \epsilon \nabla V}
  \end{equation}
  where $V\in\feszero\subset\hoz(\W)$ is a \emph{discrete variation}. Since
  $U\in\feszero$ is the discrete minimiser of the energy functional we
  certainly have that $o'(0) = 0$ and we may explicitly compute this
  quantity, the \emph{first variation},
  \begin{equation}
    \begin{split}
      o'(\epsilon) &= \int_\W \Transpose{\pd 3 L {\qp{\geovec x, U + \epsilon
              V, \nabla U + \epsilon \nabla V} }} \nabla V 
      \\
      &
      \qquad\qquad\qquad\qquad\qquad\qquad
      +
      \pd 2 L \qp{\geovec x, U + \epsilon V, \nabla U +
          \epsilon \nabla V} V \highlight{\d \geovec x} .
    \end{split}
  \end{equation}
  
  Note that since $\nabla U$ is not continuous over the skeleton of
  the triangulation $\T{}$ we have that spatial derivatives of $\pd 3
  L$ are, in general, not well defined. But $\nabla U$ is smooth over
  the interior of each element. We thus split the integral into
  elementwise contributions and integrate by parts elementwise. For
  brevity we note that the Lagrangian $L = L\qp{\geovec x, U, \nabla
    U}$ and drop the dependency.
  \begin{equation}
    \begin{split}
      o'(0)
      &=
      \sum_{K\in\T{}}\int_K 
      \Transpose{\qp{\pd 3 L}}
      \nabla V
      +
      \qp{\pd 2 L} V
      \highlight{\d \geovec x} 
      \\
      &=
      \sum_{K\in\T{}}\int_K 
      \qp{-\Div_K\qp{\pd 3 L} 
        +
        \pd2 L } V
        \highlight{\d \geovec x} 
        +
        \int_{\partial K}
        \Transpose{\qp{\pd2 L}}
        \geovec n_K V  \highlight{\d \geovec s} 
      \end{split}
    \end{equation}
where $\Div_K$ denotes an elementwise total divergence. We now use the
identity
\begin{equation}
  \sum_{K\in\T{}}
  \int_{\partial K}
  \Transpose{\qp{\pd 3 L}}
  \geovec n_K V  \highlight{\d \geovec s} 
  =
  \int_\E \jump{\pd 3 L} {V}
  \highlight{\d \geovec s} 
  ,
\end{equation}
and hence
\begin{equation}
  \begin{split}
    0 = o'(0)
    &=
    \int_\W 
    \qp{-\Div_K\qp{\pd 3 L} 
      +
      \pd2 L} V
    \highlight{\d \geovec x} 
    +
    \int_\E
    \jump{\pd3 L} V 
    \highlight{\d \geovec s} 
  \end{split}
\end{equation}
as required.
\end{Proof}

\begin{Example}[discrete Laplace's problem]
  For example the discrete Euler--Lagrange equations associated to
  Laplace's problem (Example \ref{ex:laplace}) are to find $U\in\feszero$
  such that
  \begin{equation}
    \label{eq:discrete-lap}
    \int_\W \qp{\Delta_K U + f}V \highlight{\d \geovec x} 
    +
    \int_\cE \jump{\nabla U} {V} \highlight{\d \geovec s}  = 0
    \Foreach V\in\feszero,
  \end{equation}
  where $\Delta_K$ is an elementwise Laplacian.

  Note that if $U$ is a piecewise linear function, the first term of
  (\ref{eq:discrete-lap}) is zero. Hence the discrete Laplacian can be
  completely characterised in terms of the jump of the gradient of $U$
  over the internal skeleton.
\end{Example}

\begin{Defn}[$\leb{2}(\W)$ projection operator]
  We define $\ltwoproj{}:\leb{2}(\W)\to\fes$ such that for each $w\in\leb{2}(\W)$ we have
  \begin{equation}
    \int_\W {\ltwoproj{}w}{V} \highlight{\d \geovec x}  = 
    \int_\W {w}{V} \highlight{\d \geovec x} 
    \Foreach V\in\fes.
  \end{equation}
\end{Defn}

\begin{The}[conserved quantities over $\cont{0}(\W)$--finite element spaces]
  \label{The:conserved-fe-noether}
  Let $u$ be the unique weak extrema to the minimisation problem
  \eqref{eq:variational-prob} and $U$ be its finite element
  approximation. Suppose that this problem satisfies a piecewise
  variational symmetry. Then the finite element solution satisfies the
  following
  \begin{equation}
    \begin{split}
      0
      =&
      \cN[U] 
      :=
      \int_\W
      \qp{-
        \Div_K\qp{\pd 3 L}
        +
        \pd 2 L} \ltwoproj{} \phi
      +
      L \Div{\geovec \xi}
      - 
      L \pd 1 { \geovec \xi}
      \\
      &\qquad
      -
      \Transpose{\qp{\pd 3 L}}
      \nabla U
      \Div{\geovec \xi}
      \highlight{\d \geovec x} 
      +
      \highlight{
        \int_{\cE} \jump{{\pd3 L} \ltwoproj{}\phi
        + L {\geovec \xi}}   \d \geovec s}
    +
    \int_{\pd{}\W} \Transpose{\pd 3 L}{\geovec n} \ltwoproj{}\phi \d s
    \end{split}
  \end{equation}
\end{The}

\begin{Proof}
  Recall we have an energy functional of the form 
  \begin{equation}
    \cJ[u] := \int_\W L\qp{\geovec x, u, \nabla u} \highlight{\d \geovec x}.
  \end{equation}
  We have that a finite element minimiser of this energy functional is
  continuous over the domain $\W$ but its derivative is not. \highlight{For
  simplicity and clarity, we will assume the group actions are
  separable, however, the result holds even if this is not the case.}

Using the outer variation argument from the proof of Lemma
\ref{lem:dis-EL-eqns} using $\ltwoproj \phi$ as the outer variation
\highlight{ and noting the additional boundary term arising since
  $\phi$ is not necessarily compactly supported} we need only
calculate the inner variation.
  
\highlight{
    The inner variation can be regarded as a change of variables on the
  independent variable \cite[\S 3.3]{GiaquintaHildebrandt1:1996}. In a
  similar calculation to that of the Proof of Theorem
  \ref{the:noether-weak} we let $i(\epsilon) =
  \cJ[U(\widetilde{\geovec x})]$ with $\widetilde{\geovec x} = \geovec
  x + \epsilon\geovec\xi$. We again split the integral into subdomains
  to obtain
  \begin{equation}
    \begin{split}
      0 
      &=
      \lim_{\epsilon\to 0} \frac{1}{\epsilon}\qp{\cJ[U(\widetilde{\geovec x})] - \cJ[U(\geovec x)]}
      \\
      &=
      \lim_{\epsilon\to 0}
      \frac{1}{\epsilon}
      \qp{\int_{\W_\Xi} 
        L\qp{\widetilde{\geovec x}, 
          \widetilde{U}(\widetilde{\geovec x}), 
          \widetilde{\nabla} U\qp{\widetilde{\geovec x}}
          \d\widetilde{\geovec x}
        }
        -
        \int_\W L\qp{{\geovec x}, 
          U({\geovec x}), 
          \nabla U\qp{{\geovec x}}
          \d{\geovec x}
        }
      }
      \\
      &=
      \lim_{\epsilon \to 0}
      \frac{1}{\epsilon}
      \qp{
        \int_{\W} 
        L\qp{
          \widetilde{\geovec x}
          ,
          \widetilde{U}(\widetilde{\geovec x})
          ,
          \widetilde{\nabla} U(\widetilde{\geovec x})
        }
        \frac{\d \widetilde {\geovec x}}
        {\d \geovec x} 
        \d \geovec x
        -
        \int_{\W} L\qp{\geovec x, U(\geovec x), \nabla U(\geovec x)}
        \d \geovec x 
      },
    \end{split}
  \end{equation}
  using the same coordinate transform as in the proof of Theorem
  \ref{the:noether-weak}.}
  
\highlight{
  Computing the quantities elementwise we have
  \begin{equation}
    \begin{split}
      0
      &=
      \lim_{\epsilon \to 0}
      \frac{1}{\epsilon}
      \qp{
        \sum_{K\in\T{}}
        \int_K
        L\qp{
          \widetilde{\geovec x}
          ,
          \widetilde{U}(\widetilde{\geovec x})
          ,
          \widetilde{\nabla} U(\widetilde{\geovec x})
        }
        \frac{\d \widetilde {\geovec x}}
        {\d \geovec x} 
        \d \geovec x
        -
        \int_{K} L\qp{\geovec x, U(\geovec x), \nabla U(\geovec x)}
      }
      \\
      &= 
      \sum_{K\in\T{}}
      \bigg(
      \int_K 
      - L \pd 1{\geovec \xi}
      +
      \Div{\qp{\geovec \xi}}
      \qp{L 
        -
        \Transpose{\qp{\pd 3 L}} \nabla U
      }
      \highlight{\d \geovec x}
      +
      \int_{\partial K}
      L \Transpose{\geovec \xi} \geovec n_K  
      \highlight{\d \geovec s} 
      \bigg)
      \\
        &=
        \int_\W 
        \Transpose{\qp{\nabla_K L}} \geovec \xi
        +
        L \Div{\geovec \xi}
        +
        \Transpose{\qp{\frac{\partial L}{\partial \qp{\nabla U}}}}
        \nabla U
        \Div{\geovec \xi}       \highlight{\d \geovec x} 
        +
        \int_\cE
        \jump{L \geovec \xi}      \highlight{\d \geovec s} 
        ,
      \end{split}
    \end{equation}
    as required.}
\end{Proof}

\subsection{Applications to the $p$--Laplacian}
\label{sec:p-lap}

In this section we give a numerical verification to Theorem
\ref{The:conserved-fe-noether} for a simple test problem, that of the
$p$--Laplacian
\begin{equation}
  \label{eq:p-lap}
  -\div\qp{\norm{\nabla u}^{p-2}\nabla u} = f,
\end{equation}
where we will restrict $p\in (1,\infty)$. The
$p$--Laplacian is the Euler--Lagrange equation of the following
minimisation problem: Find $u\in\sobz{1}{p}(\W)$ such that
\begin{equation}
  \label{eq:p-lap-minimise}
  \cJ_p[u] \leq \cJ_p[v] \Foreach v\in\sobz{1}{p}(\W),
\end{equation}
with the (parameterised) action functional $\cJ_p$ given by
\begin{equation}
  \cJ_p[v] := \int_\W \frac{1}{p}\norm{\nabla v}^p - fv       \highlight{\d \geovec x} .
\end{equation}

Note that for $p=2$ this problem coincides with the standard Laplace's
problem (see Example \ref{ex:laplace}). For general $p$ it is well
known that the problem is uniquely solvable.

The discrete weak formulation associated to the minimisation problem
(\ref{eq:p-lap-minimise}) is to find $U\in\fes$ such that
\begin{equation}
  \label{eq:fe-p-lap}
  \int_\W \norm{\nabla U}^{p-2}\Transpose{\qp{\nabla U}}\nabla V      \highlight{\d \geovec x} 
  =
  \int_\W f \ V       \highlight{\d \geovec x} \Foreach V\in\fes.
\end{equation}

In this test we choose $f$ such that
\begin{equation}
  \label{eq:benchmark-solution}
  u = \sin{{\pi\norm{\geovec x}^2}}
\end{equation}
solves the $p$--Laplace equation (\ref{eq:p-lap}). We have that $f$
can be written as $f = f(\norm{\geovec x})$ and hence the Lagrangian
\begin{equation}
  L(\geovec x, v, \nabla v) 
  =
  \frac{1}{p}\norm{\nabla v}^p - fv
\end{equation}
is invariant under $SO(d)$ group actions.

\begin{figure}[h]
  \caption[]{\label{fig:discrete-noether-tri-and-sol} An example of the triangulation $\T{}$ and the finite element approximation to $u = \sin{{\pi\norm{\geovec x}^2}}$, the solution of the $p$--Laplacian.}
    \subfloat[][{An example of the triangulation, here $\dim{\fes} = 12564$.}]{
      \includegraphics[scale=\figscale, width=0.47\figwidth]
                      {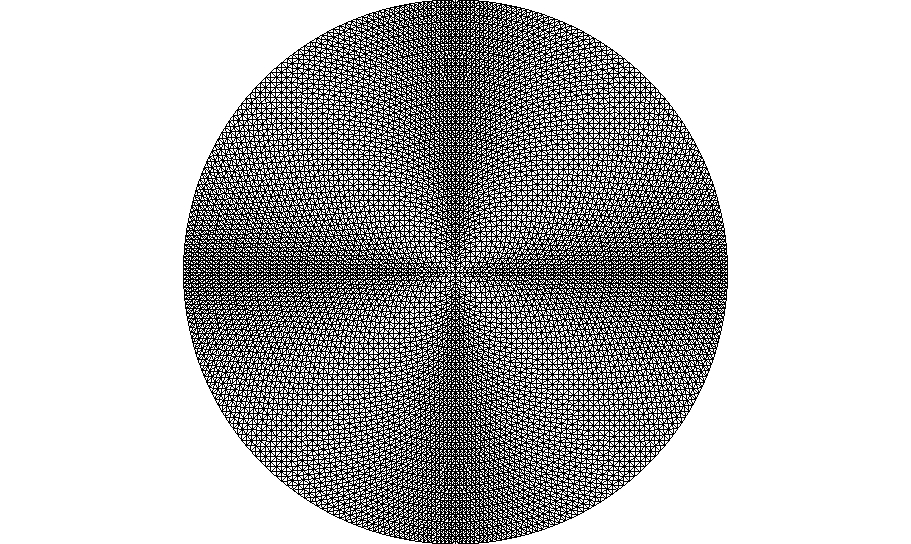}  
    }
    \hfill
    \subfloat[][{The piecewise linear finite element approximation of the solution of the $3$--Laplacian.}]{
      \includegraphics[scale=\figscale, width=0.47\figwidth]
                      {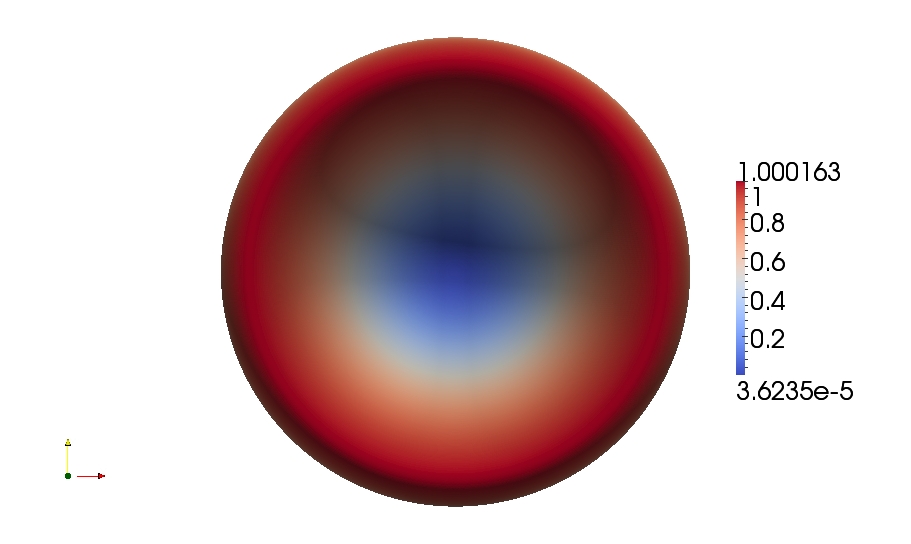}
    }
\end{figure}

We fix $d=2$ and take $\T{}$ to be a structured triangulation of $\W$,
the unit circle, as given in Figure
\ref{fig:discrete-noether-tri-and-sol}

{ It is well known \cite[c.f.]{BarrettLiu:1994} that the
  finite element approximation (\ref{eq:fe-p-lap}) is well posed and
  has optimal convergence properties. In Tables
  \ref{table:p-lap-p-3}--\ref{table:p-lap-p-5} we show errors,
  convergence rates and the values of the finite element Noether
  quantity as written in Theorem \ref{The:conserved-fe-noether} for
  various cases of $p$. The tables also study the \emph{experimental
    order of convergence} of the numerical approximation which we now define. }
\begin{Defn}[experimental order of convergence]
  \label{def:EOC}
  {Given two sequences $a(i)$ and $h(i)\downto0$, $i\integerbetween
  l{}$, we define experimental order of convergence (\EOC) to be the
  local slope of the $\log a(i)$ vs. $\log h(i)$ curve, i.e.,
  \begin{equation}
    \EOC(a,h;i):=\frac{ \log(a(i+1)/a(i)) }{ \log(h(i+1)/h(i)) }.
  \end{equation}}
\end{Defn}

\begin{table}
  \centering
  \caption{\label{table:p-lap} {In this test we computationally
      study the behaviour of the finite element conserved quantity,
      $\cN[U]$, given in Theorem \ref{The:conserved-fe-noether}. To that
      end we consider the $p$--Laplacian with various values of $p$. We fix $f$ such
      that $u$ is known and is given by
      (\ref{eq:benchmark-solution}). We compute the piecewise linear
      ($k=1$) finite element approximation given by
      (\ref{eq:fe-p-lap}).}}
  \subfloat[\label{table:p-lap-p-3} A simulation with $p=3$. \highlight{The finite element conserved quantity $\cN[U]$ is below the tolerance of the solvers.}
  ]{

\begin{tabular}{c|c|c|c|c|c} 
$\dim{\fes}$ & $\Norm{u - U^N}_{\leb{p}(\W)}$ & EOC & $\Norm{u-U^N}_{\sob{1}{p}(\W)}$ & EOC & $\cN[U]$ \\ 
\hline \\ 
13 & 1.12725161 & 0.000 & 3.63060492 & 0.000 &2.353673e-14 \\ 
 41 & 0.70442091 & 0.678 & 2.98329491 & 0.283 &1.776357e-14 \\ 
 145 & 0.15390246 & 2.194 & 1.56695787 & 0.929 &3.108624e-15 \\ 
 585 & 0.03539738 & 2.120 & 0.73020113 & 1.102 &1.526557e-13 \\ 
 2805 & 0.00618342 & 2.517 & 0.31109110 & 1.231& 2.456799e-12 \\ 
 14293 & 0.00110803 & 2.480 & 0.13244566 & 1.232&  2.600439e-12 \\ 
 73401 & 0.00022205 & 2.319 & 0.05750230 & 1.204 & 6.134856e-13 \\ 
 384769 & 0.00004666 & 2.250 & 0.02645758 & 1.120 & 1.581680e-13 \\ 
 \end{tabular}
  }
  \\
  \subfloat[\label{table:p-lap-p-4} A simulation with $p=4$. \highlight{The finite element conserved quantity $\cN[U]$ is below the tolerance of the solvers.}
  ]{
\begin{tabular}{c|c|c|c|c|c} 
$\dim{\fes}$ & $\Norm{u - U^N}_{\leb{p}(\W)}$ & EOC & $\Norm{u-U^N}_{\sob{1}{p}(\W)}$ & EOC & $\cN[U]$\\ 
\hline \\ 
13 & 1.50645187 & 0.000 & 3.89649998 & 0.000 &4.002354e-14 \\ 
 41 & 0.95656811 & 0.655 & 3.32660211 & 0.228 &1.776357e-14 \\ 
 145 & 0.18567964 & 2.365 & 1.69441051 & 0.973 &3.730349e-14 \\ 
 585 & 0.04346809 & 2.095  &0.77624620 & 1.126 &1.314726e-12 \\ 
 2805 & 0.00787713 & 2.464 & 0.33055071 & 1.232 &2.027223e-11 \\ 
 14293 & 0.00141755 & 2.474 & 0.13947531 & 1.245 &2.106114e-11 \\ 
 73401 & 0.00028203 & 2.329 & 0.06028734 & 1.210 &4.074260e-12 \\ 
 384769 & 0.00005934 & 2.249 &0.02756628 & 1.129 &1.291921e-12 \\ 
 \end{tabular}
  }
  \\
  \subfloat[\label{table:p-lap-p-5} A simulation with $p=5$. \highlight{The finite element conserved quantity $\cN[U]$ is below the tolerance of the solvers.}
  ]{\begin{tabular}{c|c|c|c|c|c} 
$\dim{\fes}$ & $\Norm{u - U^N}_{\leb{p}(\W)}$ & EOC & $\Norm{u-U^N}_{\sob{1}{p}(\W)}$ & EOC & $\cN[U]$\\ 
\hline \\ 
13 & 1.84301273 & 0.000 &4.16722611 & 0.000 &7.812362e-13 \\ 
 41 & 1.14195690 & 0.691& 3.69726521 & 0.173 &8.526513e-14 \\ 
 145 & 0.21243009 & 2.426& 1.85724530 & 0.993 &2.984279e-13 \\ 
 585 & 0.05034349 & 2.077 &0.83769960 & 1.149 &1.120704e-11 \\ 
 2805 & 0.00932480 & 2.433 &0.35992898 & 1.219 &1.558869e-10 \\ 
 14293 & 0.00170255 & 2.453 &0.15093728 & 1.254 &1.587779e-10 \\ 
 73401 & 0.00033701 & 2.337 &0.06531226 & 1.209 &2.558158e-11 \\ 
 384769 & 0.00007080 & 2.251 &0.02955021 & 1.144 &1.013448e-11 \\ 
 \end{tabular}
  }
\end{table}

\begin{Rem}[\highlight{numerical conservation}]
  \highlight{ In the numerical experiments conducted in Tables
    \ref{table:p-lap-p-3}--\ref{table:p-lap-p-5} we formulated
    (\ref{eq:fe-p-lap}) as a system of nonlinear equations, the
    solution to this is then approximated by a Newton method with
    tolerance set at $10^{-10}$. At each Newton step the solution to
    the linear system of equations is approximated using a stabilised
    conjugate gradient iterative solver with an algebraic multigrid
    preconditioner, also set at a tolerance of $10^{-10}$. Since the
    solvers themselves only generate approximations to the numerical
    variational problem, the notion of conservation is only true up to
    a certain tolerance. In this case, the quantity will be conserved
    up to the tolerance of the solvers, $10^{-10}$.}
\end{Rem}

\subsection{{Mimetic methods weakly enforce discrete conservation laws which are derived from trivial Lie group actions}}

The mimetic finite element framework consists of reformulating the
Euler--Lagrange equations as a system of first order PDEs. Consider
our prototypical example for illustrative purposes. Poisson's problem,
\begin{equation}
  \Delta u = 0,
\end{equation}
is the Euler--Lagrange equation of the minimisation problem
\begin{equation}
  \cJ[u] = \int_\W \frac{1}{2}\norm{\nabla u}^2       \highlight{\d \geovec x} \to \text{ min.}
\end{equation}
It can be written in \emph{mixed form} by introducing an auxiliary
variable $\geovec p$ to represent the gradient and rewriting Poisson's
problem to seek $\qp{u, \geovec p}$ such that
\begin{gather}
  \div{\geovec p} = 0
  \\
  \geovec p = \nabla u.
\end{gather}
These are then the Euler--Lagrange equations of the saddle point problem
\begin{equation}
  \cK[u, \geovec p] 
  := 
  \int_\W \frac{1}{2} \norm{\geovec p}^2
  +
  u\qp{\div{\geovec p}}       \highlight{\d \geovec x}.
\end{equation}
The correct function space setting is to seek $u\in\leb{2}(\W)$ and
$\geovec p \in\sobh{\text{div}}(\W) := \ensemble{\geovec \Psi
}{\div{\geovec \Psi} \in \leb{2}(\W)}$. A conformal approximation of
this problem can be sought using the Raviart--Thomas and piecewise
constant finite element pair \cite{RaviartThomas:1977}, for example. A
sufficient condition for the construction of a conformal finite
element space of $\sobh{\text{div}}(\W)$ is that the jumps of the
discrete functions vanish over the skeleton of the domain
\cite[c.f.]{BrezziFortin:1991}. 

Recall Remark \ref{rem:trivial-actions} concerned itself with the
trivial Lie group action of translation in the dependent variable. For
our model problem we have that $\nabla u$ is a conservation law. The
mimetic scheme weakly enforces this conservation law.

\section{Conservative properties of Lagrangian FEs for strong solutions}
\label{sec:numerics-lag}

In this section we present results concerning the approximability of
the strong continuous conservation laws arising from Theorem
\ref{the:noether}. We examine numerically the behaviour of the
Lagrangian finite element method. In this sense we wish to measure the
quantity $\Div\qp{\cC[U]}$ and evaluate how far it deviates from
zero. For clarity of exposition, we will assume henceforth that the
continuous minimisation problem takes the form: Find $u\in\hoz(\W)$ such
that
\begin{equation}
  \cJ[u] = \inf_{v\in\hoz(\W)} \cJ[v].
\end{equation}

\begin{The}[Bound on the finite element approximation of Noether's laws.]
  \label{the:pre-convergence}
  Let $u\in\sobh{2}(\W)\cap \hoz(\W)$ be a strong extrema to the
  variational problem (\ref{eq:variational-prob}) (and hence a strong
  solution to the Euler--Lagrange equations
  (\ref{eq:Euler-Lagrange-equations})). Suppose we have that Theorem
  \ref{the:noether} holds under a variational symmetry group with
  infinitesimals $\geovec \xi$ and $\phi$.

  In addition assume that we have 
\begin{gather}
  \label{eq:ass-3}
  L(\geovec x, u, \nabla u) \in \leb{\infty}(\W) 
  \Foreach u\in\hoz(\W)
  \\
  \label{eq:ass-4}
  \pd 3 L(\geovec x,u,\nabla u)
  \in \qb{\leb{\infty}(\W)}^d
  \Foreach u\in\hoz(\W).
\end{gather}
Then if $U\in\fes$ is the finite element approximation to $u$ there
exists a constant $C$ such that
\begin{equation}
  \label{eq:the-bound}
  \begin{split}
  \Norm{\Div{\cC[U]}}_{\sobh{-1}(\W)} 
  &\leq
  C
  \Bigg(
    \Norm{L(\geovec x, U, \nabla U)
      -
      L(\geovec x, u, \nabla u)
    }_{\leb{2}(\W)}
    +
    \Norm{
    \nabla U - \nabla u
    }_{\leb{2}(\W)}
    \\
    &
    \qquad
    +
    \Norm{\pd 3 L(\geovec x,U,\nabla U)
      -
      \pd 3 L(\geovec x,u,\nabla u)
    }_{\leb{2}(\W)}
    \\
    &
    \qquad
    +
    \Norm{\phi(\geovec x, U) - \phi(\geovec x, u)}_{\leb{2}(\W)}
    +
    \Norm{\geovec \xi(\geovec x, U) - \geovec \xi(\geovec x, u)}_{\leb{2}(\W)}
  \Bigg).
  \end{split}
  \end{equation}
\end{The}
\begin{Proof}
  We begin by noting that since $u$ is a strong extremal, Theorem
  \ref{the:noether} holds and we have that $\Div\qp{\cC[u]} = 0$. Hence
  \begin{equation}
    \Norm{\Div{\cC[U]}}_{\sobh{-1}(\W)}
    =
    \Norm{\Div\qp{\cC[U] - \cC[u]}}_{\sobh{-1}(\W)}.
  \end{equation}
  Now, we may use the fact that for a generic $\varphi\in\hoz(\W)$
  \begin{equation}
    \begin{split}
      \duality{\Div\qp{\cC[U] - \cC[u]}}{\varphi}
      &=
      -\ltwop{{\cC[U] - \cC[u]}}{\nabla \varphi}
      \\
      &\leq 
      \Norm{{\cC[U] - \cC[u]}}_{\leb{2}(\W)}
      \Norm{\nabla\varphi}_{\leb{2}(\W)}.
    \end{split}
  \end{equation}
  Since $\varphi$ was generic we may divide through by
  $\Norm{\nabla\varphi}$ and take the supremum over $\varphi$ then by
  the definition of the $\sobh{-1}(\W)$ norm
  \begin{equation}
    \Norm{\Div{\cC[U]}}_{\sobh{-1}(\W)}
    \leq
    \Norm{\cC[U] - \cC[u]}_{\leb{2}(\W)}.
  \end{equation}
  By the definition of the Noether quantity $\cC$ from (\ref{eq:noether-functional}) we have that 
  \begin{equation}
    \begin{split}
      \Norm{\Div{\cC[U]}}_{\sobh{-1}(\W)}      
      &\leq
      \Norm{
        \phi(\geovec x, U)
        \pd 3 L(\geovec x,U,\nabla U)
        -
        \phi(\geovec x, u)
        \pd 3 L(\geovec x,u,\nabla u)
      }_{\leb{2}(\W)}
      \\
      &+
      \Norm{
        \Transpose{\qp{\geovec \xi(\geovec x, U)}}
        \nabla U
        \pd3 L(\geovec x,U,\nabla U)
        -
        \Transpose{\qp{\geovec \xi(\geovec x, u)}}
        \nabla u
        \pd 3 L(\geovec x,u,\nabla u)
      }_{\leb{2}(\W)}
      \\
      &+
      \Norm{
        L(\geovec x,U,\nabla U)
        \geovec \xi(\geovec x, U)
        -
        L(\geovec x,u,\nabla u)
        \geovec \xi(\geovec x, u)
      }_{\leb{2}(\W)}
      \\
      &=:
      \cI_1 + \cI_2 + \cI_3
    \end{split}
  \end{equation}
  where for clarity we have written the dependencies explicitly. Now
  for each of the $\cI_i$ we add and subtract appropriate quantities
  and make use of the triangle inequality. We thus have the following
  bounds:
  \begin{equation}
    \begin{split}
      \cI_1
      &\leq 
      \Norm{
        \phi(\geovec x, U)
        \qp{
          \pd3 L(\geovec x,U,\nabla U)
          -
          \pd3 L(\geovec x,u,\nabla u)
        }
      }_{\leb{2}(\W)}
      \\
      &
      \qquad 
      +
      \Norm{
        \qp{
          \phi(\geovec x, U) - \phi(\geovec x, u)
        }
        \pd3 L(\geovec x,u,\nabla u)
      }_{\leb{2}(\W)}
    \end{split}
  \end{equation}
  \begin{equation}
    \begin{split}
      \cI_2
      &\leq
      \Norm{
        \qp{
          \geovec \xi(\geovec x, u)
          -
          \geovec \xi(\geovec x, U)
        }
        \Transpose{\qp{
            \nabla u
          }}
        \pd3 L(\geovec x,u,\nabla u)
      }_{\leb{2}(\W)}
      \\
      &
      \qquad
      +
      \Norm{
        \Transpose{
          \qp{\geovec \xi(\geovec x, U)}
        }
        \qp{
          \nabla u - \nabla U
          }
        \pd3 L(\geovec x,u,\nabla u)
      }_{\leb{2}(\W)}
      \\
      &
      \qquad
      +
      \Norm{
       \Transpose{
         \qp{\geovec \xi(\geovec x, U)}
       }
       \nabla U
       \qp{
         \pd3 L(\geovec x,U,\nabla U)
         -
         \pd3 L(\geovec x,u,\nabla u)
       }
      }_{\leb{2}(\W)}
    \end{split}
  \end{equation}
  and
  \begin{equation}
    \begin{split}
      \cI_3
      &\leq
      \Norm{
        \qp{
          L(\geovec x,U,\nabla U)
          -
          L(\geovec x,u,\nabla u)
          }
        \geovec \xi(\geovec x, U)  
      }_{\leb{2}(\W)}
      \\
      &
      \qquad
      + 
      \Norm{
        L(\geovec x,u,\nabla u)
        \qp{
          \geovec \xi(\geovec x, U) 
          -
          \geovec \xi(\geovec x, u)  
        }
      }_{\leb{2}(\W)}.
    \end{split}
  \end{equation}
  
Under assumptions (\ref{eq:ass-3})--~(\ref{eq:ass-4}) and using the
fact that $\phi$ and $\geovec \xi$ are infinitesimals of a smooth Lie
group action we have that
\begin{equation}
  \begin{split}
    \cI_1 
    &\leq
    C \qp{ \Norm{ \pd 3 L(\geovec x,U,\nabla
        U)
        - \pd3 L(\geovec
        x,u,\nabla u) }_{\leb{2}(\W)} 
      +
      \Norm{ \phi(\geovec x, U) - \phi(\geovec x, u) }_{\leb{2}(\W)} }
    \\
    \cI_2
    &\leq
    C\bigg(\Norm{ \geovec \xi(\geovec x, u) - \geovec
        \xi(\geovec x, U) }_{\leb{2}(\W)} + \Norm{ \nabla u - \nabla U
      }_{\leb{2}(\W)}
      \\
      &\qquad\qquad\qquad\qquad\qquad\qquad + \Norm{ \pd3 L(\geovec x,U,\nabla
        U)
        -
        \pd3 L(\geovec
        x,u,\nabla u)
      }_{\leb{2}(\W)}\bigg)
      \\
      \cI_3 &\leq C \qp{ \Norm{ L(\geovec x,U,\nabla U) - L(\geovec
        x,u,\nabla u) }_{\leb{2}(\W)} + \Norm{ \geovec \xi(\geovec x,
        u) - \geovec \xi(\geovec x, U) }_{\leb{2}(\W)} }.
  \end{split}
\end{equation}
Taking the sum of the $\cI_i$ gives the desired result.
\end{Proof}

\begin{Cor}
  \label{cor:under-lipscitz-ass}
{
    Let the conditions of Theorem \ref{the:pre-convergence} hold under
  the same variational symmetry group with infinitesimals $\geovec
  \xi$ and $\phi$. In addition assume that the Lagrangian is
  sufficiently smooth such that both $L$ and $\frac{\partial
    L}{\partial(\nabla u)}$ are (locally) Lipschitz with respect to
  the second and third variable then the bound (\ref{eq:the-bound})
  can be simplified to}
\begin{equation}
  \begin{split}
  \Norm{\Div{\cC[U]}}_{\sobh{-1}(\W)} 
  &\leq
  C
  \Bigg(
    \Norm{
      U - u
    }_{\leb{2}(\W)}
    +
    \Norm{
    \nabla U - \nabla u
    }_{\leb{2}(\W)}
    \\
    &
    \qquad
    +
    \Norm{\phi(\geovec x, U) - \phi(\geovec x, u)}_{\leb{2}(\W)}
    +
    \Norm{\geovec \xi(\geovec x, U) - \geovec \xi(\geovec x, u)}_{\leb{2}(\W)}
  \Bigg).
  \end{split}
  \end{equation}
\end{Cor}

\begin{Rem}
  \highlight{ The results of Theorem \ref{the:pre-convergence} are not
    just applicable to the finite element solution, but to any
    function. Indeed, the result is actually a property of the
    conservation law $\cC[\cdot]$ rather than the approximation $U$.
  }
\end{Rem}

\begin{Rem}[relating to the $p$--Laplacian]
  \label{rem:convergence}
  We may relate Theorem \ref{the:pre-convergence} and Corollary
  \ref{cor:under-lipscitz-ass} to the $p$--Laplacian studied in
  \S\ref{sec:p-lap}. We were considering that $L$ was invariant under
  rotations in the independent variable. In that case we have that
  $\phi \equiv 0$ and that
  \begin{equation}
    \Norm{\geovec \xi(\geovec x, U) - \geovec \xi(\geovec x, u)}
    =
    0
  \end{equation}
  hence we may infer that 
  \begin{equation}
    \Norm{\Div{\cC[U]}}_{\sobh{-1}(\W)} 
    \leq
    C
    \bigg(
    \Norm{
      U - u
    }_{\leb{2}(\W)}
    +
    \Norm{
      \nabla U - \nabla u
    }_{\leb{2}(\W)}
    \bigg).
  \end{equation}
  The leading term here is $\Norm{\nabla U - \nabla u}_{\leb{2}(\W)}$
  which is well known to be = $\Oh(h^k)$.
\end{Rem}

\begin{Lem}[a computable upper bound for the conservation law]
  \label{lem:upper-bound}
  Let $u\in\sobh{2}(\W)\cap \hoz(\W)$ be a strong extrema to the
  variational problem (\ref{eq:variational-prob}). Suppose that
  Theorem \ref{the:noether} holds and that $U$ is the finite element
  approximation to $u$. Then there exists a constant $C$ dependent on
  the shape regularity of $\T{}$ such that
  \begin{equation}
    \Norm{\Div{\cC[U]}}_{\sobh{-1}(\W)} 
    \leq
    E\qp{U,f} := 
    C \qp{
    \sum_{K\in\T{}} \Norm{\Div{\cC[U]}}_{\leb{2}(K)}
    +
    \sum_{e\in\E} \Norm{\jump{\cC[U]}}_{\leb{2}(e)}
    }.
  \end{equation}
\end{Lem}
\begin{Proof}
  The proof is a standard aposteriori argument where the quantity of
  interest is split into regular and singular parts. Recall for the
  model elliptic problem $\cC[U] = \cC(\geovec x, U, \nabla U)$ and
  hence $\Div{\cC[U]} \notin \leb{2}(\W)$. Let
  $\duality{\cdot}{\cdot}$ denote the $\sobh{-1}(\W)$ -- $\hoz(\W)$
  duality pairing then for any $\varphi\in\hoz(\W)$ it holds that
  \begin{equation}
    \begin{split}
      \duality{\Div{\cC[U]}}{\varphi}
      &=
      \sum_{K\in\T{}}
      - \int_K \Transpose{\cC[U]}\nabla \varphi      \highlight{\d \geovec x} 
      \\
      &=
      \sum_{K\in\T{}}
      \int_K \Div{\cC[U]}\varphi       \highlight{\d \geovec x} 
      -
      \int_{\partial K} \Transpose{\cC[U]}\geovec n_K \varphi      \highlight{\d \geovec x}.
    \end{split}
  \end{equation}
  Note that for generic $\geovec p \in \qb{\leb{2}(\W)}^d$ and $v\in\sobh{1}(\W)$ that
  \begin{equation}
    \sum_{K\in\T{}} 
    \int_{\partial K} 
    \Transpose{\geovec p} \geovec n_K v = \int_\E \jump{\geovec p} v       \highlight{\d \geovec x} 
    +
    \int_{\partial \W} \Transpose{\geovec p} \geovec n_K v       \highlight{\d \geovec s} .
  \end{equation}
  Hence
  \begin{equation}
    \duality{\Div{\cC[U]}}{\varphi} 
    =
    \sum_{K\in\T{}} 
    \int_K \Div{\cC[U]}\varphi       \highlight{\d \geovec x} 
    -
    \sum_{e\in\E}
    \int_e\jump{\cC[U]} \varphi      \highlight{\d \geovec s} 
    -
    \sum_{e\in\partial\W}
    \int_e\Transpose{\cC[U]}\geovec n \varphi       \highlight{\d \geovec s} .
  \end{equation}
  Applying a Cauchy--Schwarz inequality followed by a Poincar\'e
  inequality together with a trace inequality yields
  \begin{equation}
    \begin{split}
      \duality{\Div{\cC[U]}}{\varphi} 
      &\leq
      C\sum_{K\in\T{}}
      \Norm{\Div{\cC[U]}}_{\leb{2}(K)} \Norm{\varphi}_{\leb{2}(K)} 
      +
      \sum_{e\in\E}
      \Norm{\jump{\cC[U]}}_{\leb{2}(e)}\Norm{\varphi}_{\leb{2}(e)}
      \\
      &\leq
      C\Norm{\nabla\varphi}_{\leb{2}(\W)}
      \qp{ \sum_{K\in\T{}}
      \Norm{\Div{\cC[U]}}_{\leb{2}(K)} 
      +
      \sum_{e\in\E}
      \Norm{\jump{\cC[U]}}_{\leb{2}(e)}}.
    \end{split}
  \end{equation}
  Noting $\varphi$ was a generic function in $\hoz(\W)$, dividing through by
  $\Norm{\nabla\varphi}_{\leb{2}(\W)}$ and taking the supremum over $\varphi$ yields the
  desired result.
\end{Proof}

\begin{Rem}[]
  It is worth noting that since the quantity $\Div \cC[U]$ is not
  orthogonal to $\fes$ there are no powers of $h$ appearing in the
  estimate given in Lemma \ref{lem:upper-bound}. What is interesting
  is that the estimate still converges to zero at the same rate as
  that of the residual of the problem \highlight{in an aposteriori
    sense \cite[c.f.]{Ainsworth:2000}}.
\end{Rem}

\subsection{Numerical experiments}
\label{sec:benchmark}

\begin{figure}[h!]
  \caption[]{\label{fig:benchmark-experiements} In this experiment we
    consider the $2$--Laplacian. We fix $f$ such that $u =
    \exp\qp{-10\norm{\geovec x}^2}$. We solve the discrete problem on
    concurrently refined meshes and compute the $\leb{2}(\W)$--error,
    the $\smash{\hoz(\W)}$--error and the computable estimate given in Lemma
    \ref{lem:upper-bound}. Notice that in each of the examples
    $E(U,f)$ converges like $\Oh\qp{h^k}$ \highlight{as predicted in
      Remark \ref{rem:convergence}}.}
  \begin{center}
    \subfloat[][{Initial mesh}]{
      \includegraphics[scale=\figscale, width=0.47\figwidth]
                      {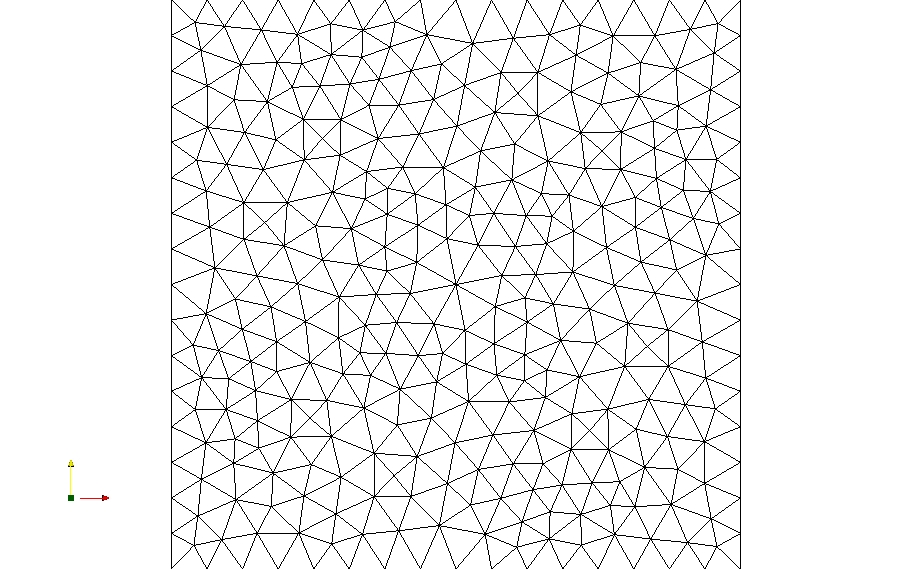}
    }
    \subfloat[][{In this example we choose $k=1$. Note when
      $\dim{\fes} = 295681, E(U,f) \approx 0.08$.}]{
      \includegraphics[scale=\figscale, width=0.47\figwidth]
                      {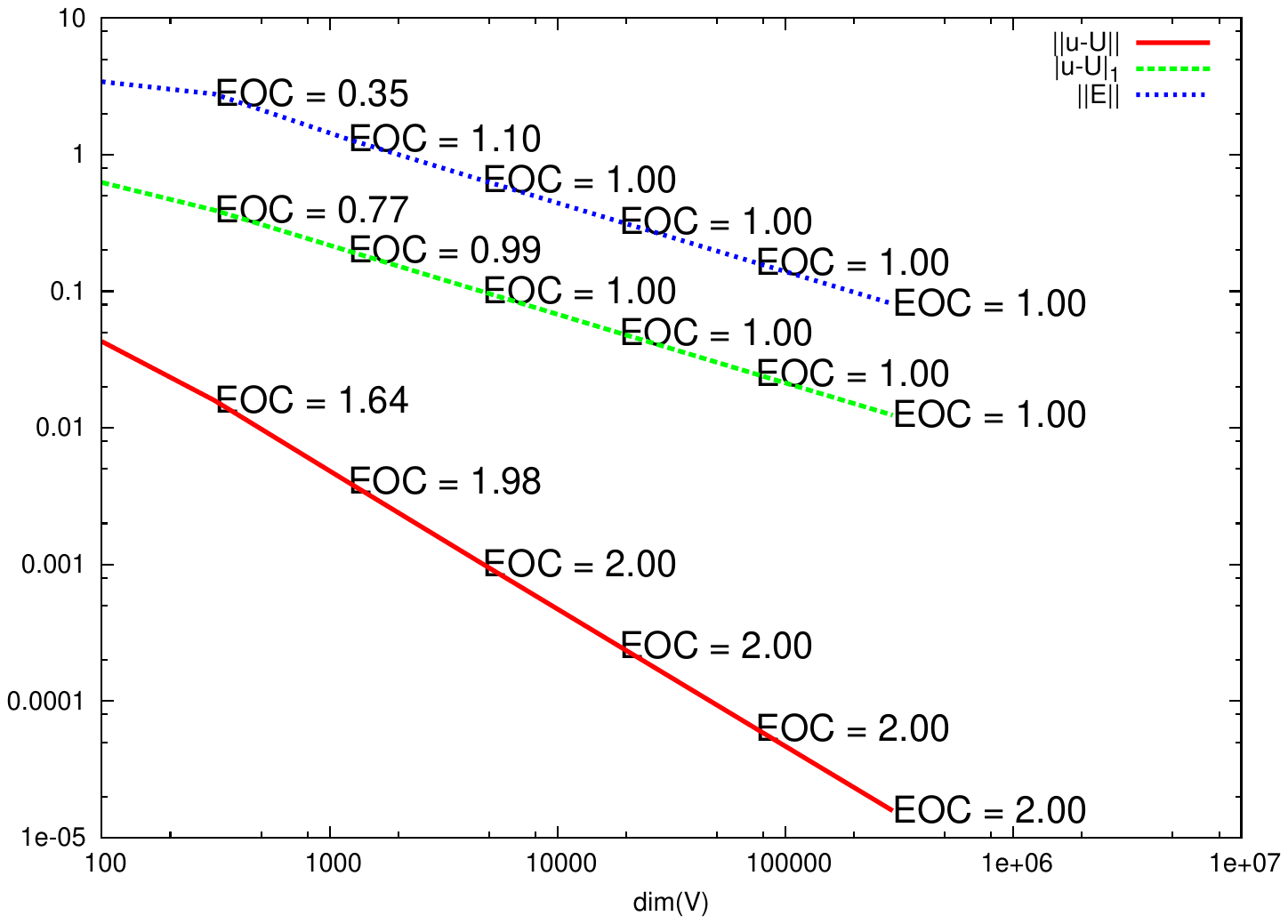}
    }\\
    \subfloat[][{In this example we choose $k=2$. Note when
      $\dim{\fes} = 295681, E(U,f) \approx 0.0003$.}]{
      \includegraphics[scale=\figscale, width=0.47\figwidth]
                      {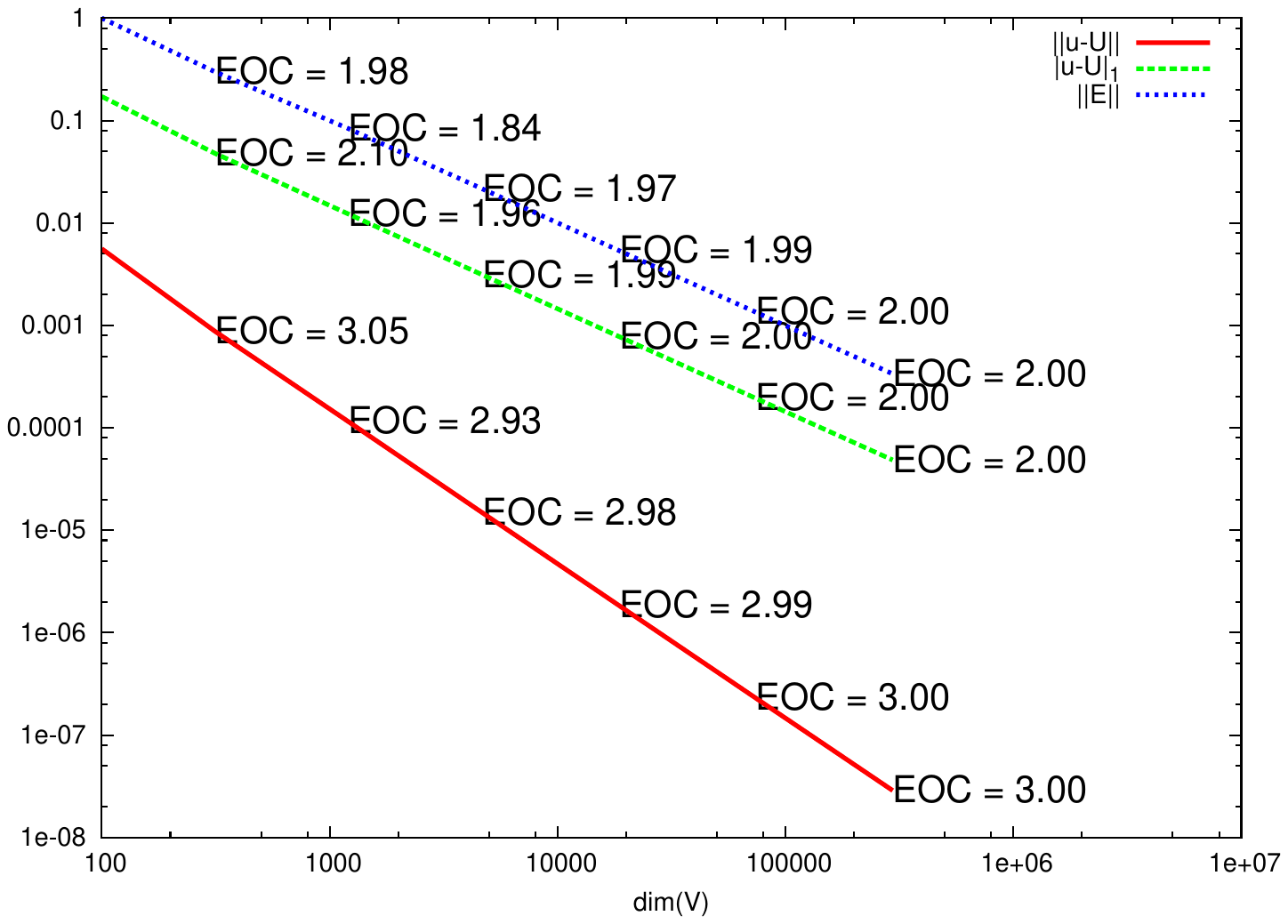}
                    }
    \subfloat[][{In this example we choose $k=3$. Note when
      $\dim{\fes} = 295681, E(U,f) \approx 0.000001$.}]{
      \includegraphics[scale=\figscale, width=0.47\figwidth]
                      {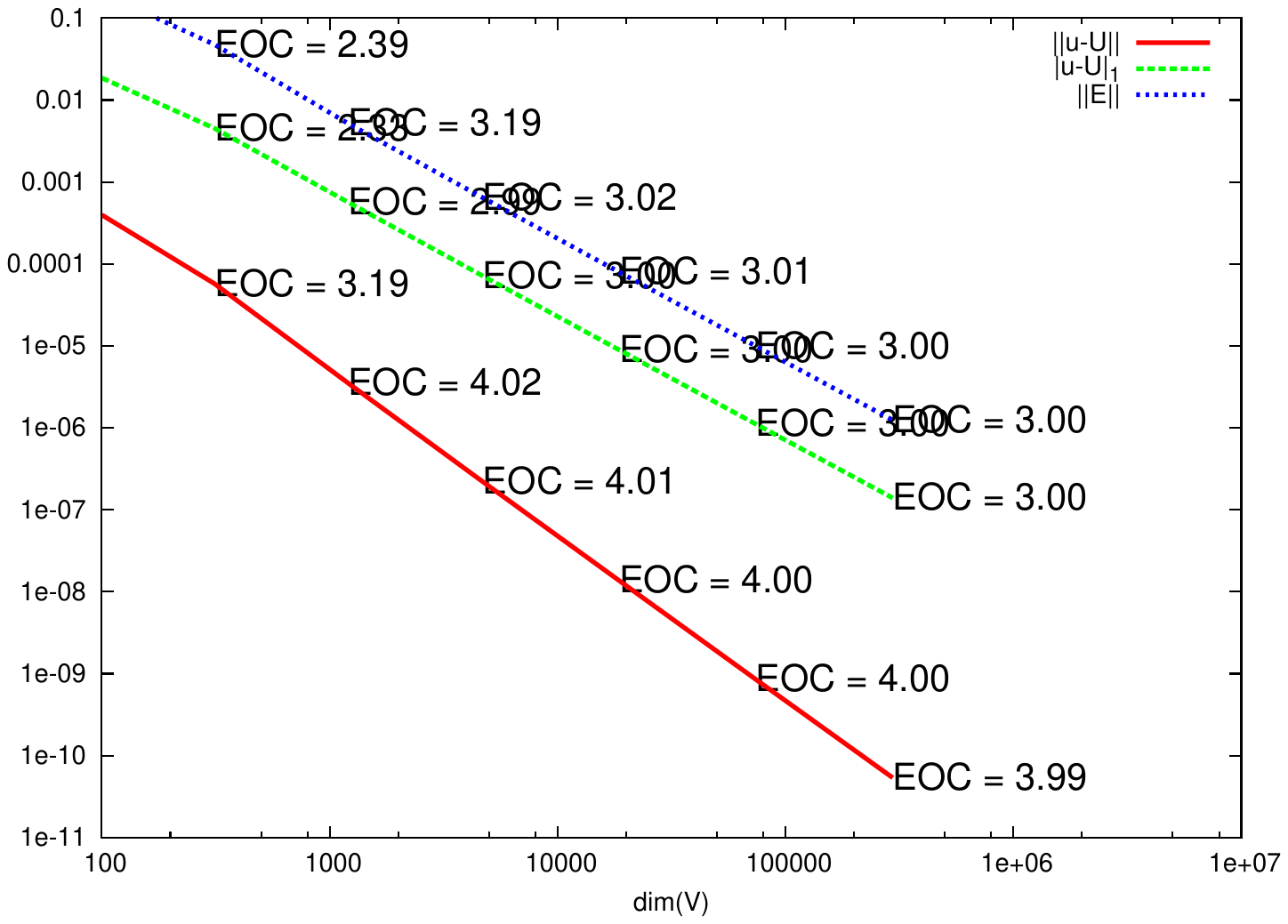}
    }
  \end{center}
\end{figure}

In Figure \ref{fig:benchmark-experiements} we show numerically that
the estimate given in Lemma \ref{lem:upper-bound} converges at the
same rate as the apriori bound given in Remark
\ref{rem:convergence}. All numerical experiments are conducted for
simplicity on the $2$--Laplacian taking $\W = [-1,1]^2$ which is
discretised using an unstructured triangulation.

\subsection{Adaptive methods}

\begin{figure}[h]
  \caption[]{\label{fig:adaptive-experiements} In this experiment we
    fix $f$ such that $u$ is as in Figure
    \ref{fig:benchmark-experiements}. We solve the FE problem on
    concurrently refined meshes and compute the $\leb{2}(\W)$--error,
    the $\smash{\hoz(\W)}$--error and the estimate of Noether's
    conservation law, $E\qp{U,f}$, from Lemma
    \ref{lem:upper-bound}. We construct an h--adaptive approximation
    to the problem with the aim of minimising the \highlight{estimate
      given in Lemma \ref{lem:upper-bound} and hence the} violation of the Noether
    quantity. Note that the estimator decreases far quicker than
    $\Oh(N^{-k/2})$.}
  \subfloat[][{Adaptive solution}]{ \includegraphics[scale=\figscale,
    width=0.47\figwidth] {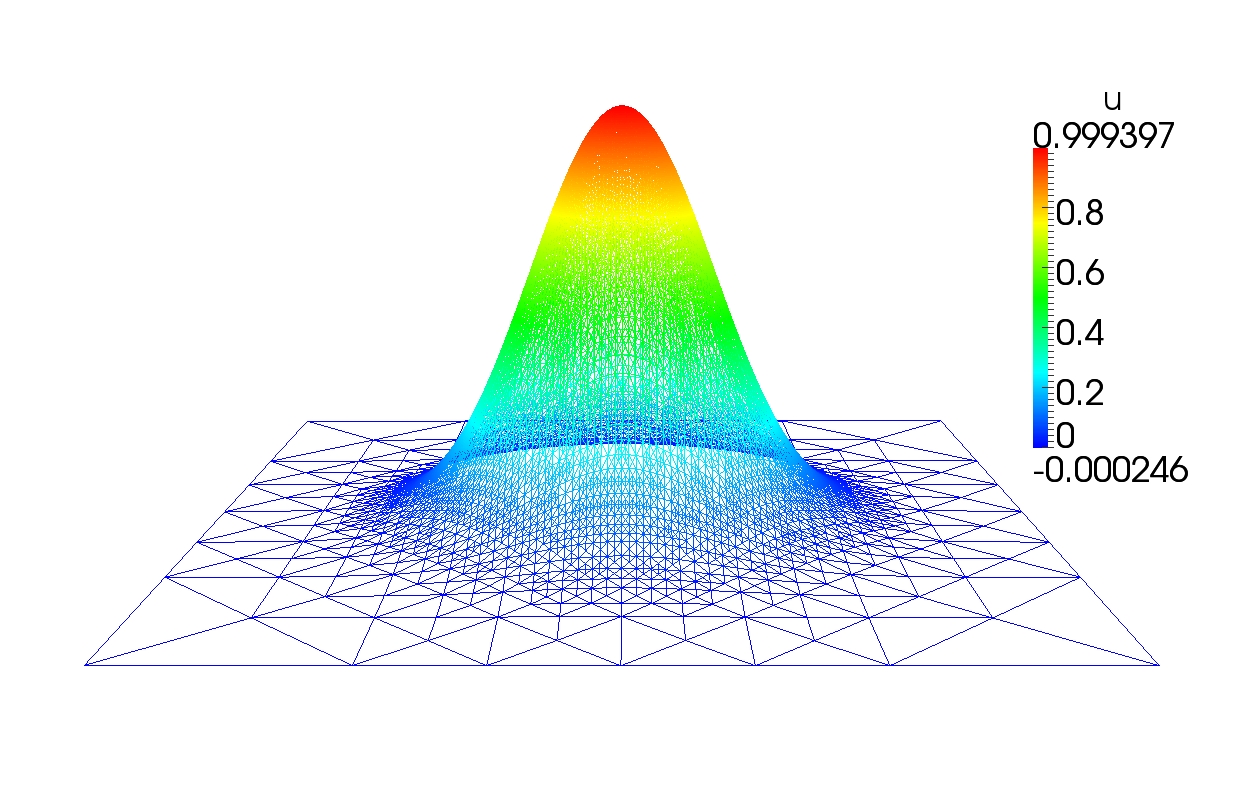} }
  \hfill \subfloat[][{In this example we choose $k=1$. Note when
    $\dim{\fes} = 40020, E(U,f) \approx 0.002$.}]{
    \includegraphics[scale=\figscale, width=0.47\figwidth]
                    {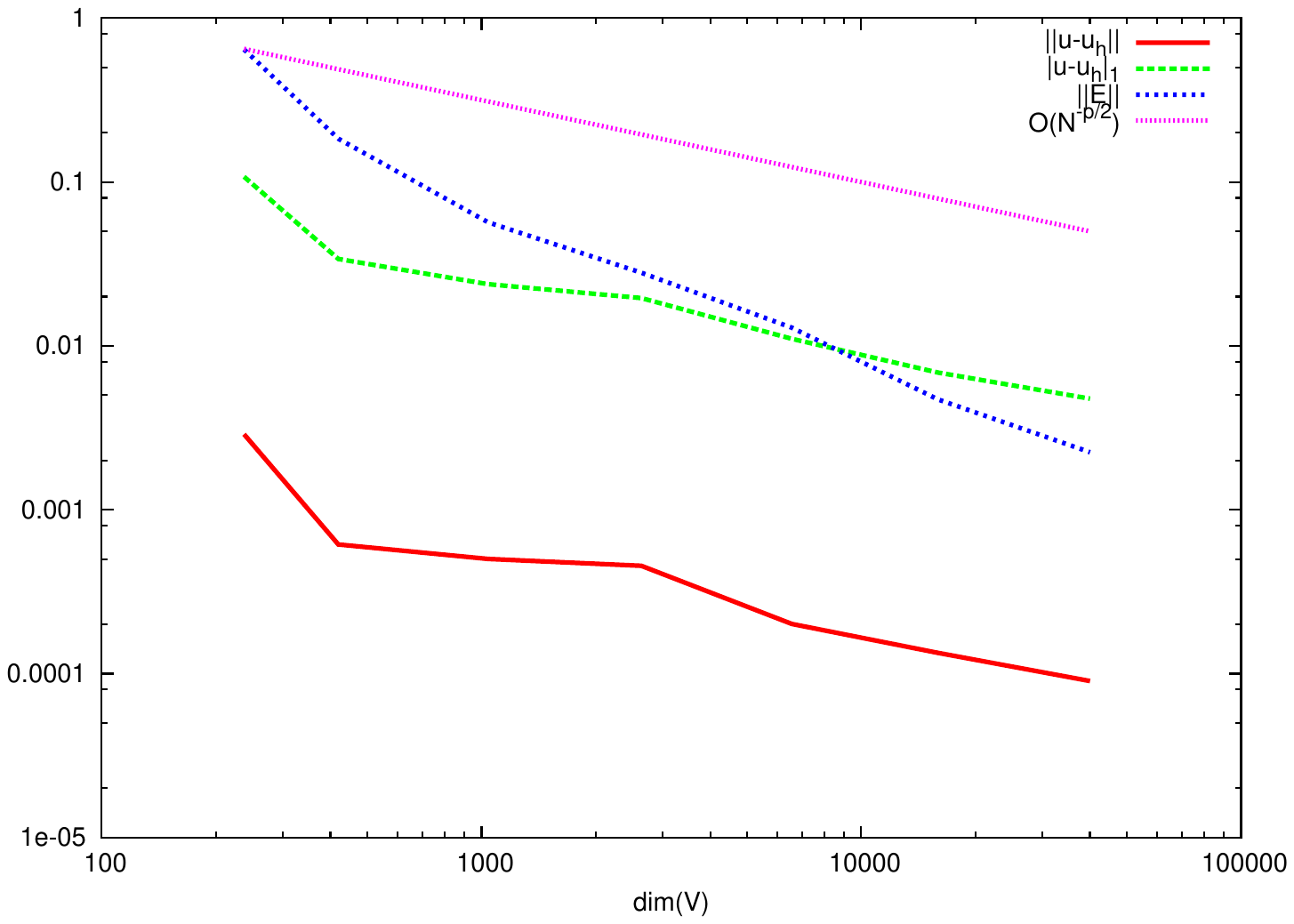} }
\end{figure}

We conclude our numerical experiments in Figure
\ref{fig:adaptive-experiements} by using the observed convergence of
the computable estimator $E(U,f)$ to construct an adaptive scheme
aimed at minimising $E(U,f)$ (and hence
$\Norm{\Div{\cC[U]}}_{-1}$). The adaptive algorithm we make use of is
of standard type (SOLVE $\to$ ESTIMATE $\to$ MARK $\to$ REFINE
\cite[c.f.]{alberta}) utilising the maximum strategy marking and
newest vertex bisection refinement.

\section{Conclusions and outlook}

In this work we have given a concise statement of Noether's first
theorem, applied to a set of model problems.

We have proved a Noether type theorem for a specific class of weak
extremum to a model variational problem, that is, those that are
$\sob{1}{\infty}(\W)$ with finitely many jump discontinuities.  We
have in addition proved an equivalent discrete theorem for the finite
element approximation of the problem. We write the exact conserved
quantity for the discrete scheme in the same spirit as
\cite{HydonMansfield:2004}. We have demonstrated that the Lagrangian
finite elements enjoy the property of asymptotically conserving the
strong Noether conservation laws when approximating strong solutions
of certain classes of variational problem and Lie group
action. \highlight{Note that although there is no requirement that the
  finite element mesh have the underlying symmetry, the weak conserved
  quantities are independent of the mesh.}

In addition we have studied the exact discrete conserved quantities
numerically. \highlight{These are conserved irrespective of whether
  the underlying symmetry is inbuilt into underlying function space.}
We have also constructed a geometric based adaptive scheme to conserve
the approximate continuous quantities up to user specified tolerance.
\highlight{ This means that upon each adaptive step there is a
  discrete conserved quantity which can be taken as close to the
  continuous counterpart as the user specifies.}

\bibliographystyle{alpha} 
\bibliography{tristansbib,tristanswritings}

\begin{thebibliography}{CMKO11}

\bibitem[AFW10]{ArnoldFalkWinther:2010}
Douglas~N. Arnold, Richard~S. Falk, and Ragnar Winther.
\newblock Finite element exterior calculus: from {H}odge theory to numerical
  stability.
\newblock {\em Bull. Amer. Math. Soc. (N.S.)}, 47(2):281--354, 2010.

\bibitem[Ang12]{Angier:2012}
Natalie Angier.
\newblock The mighty mathematician you've never heard of.
\newblock {\em The New York Times}, 2012.

\bibitem[BBL09]{BrezziBuffaLipnikov:2009}
Franco Brezzi, Annalisa Buffa, and Konstantin Lipnikov.
\newblock Mimetic finite differences for elliptic problems.
\newblock {\em M2AN Math. Model. Numer. Anal.}, 43(2):277--295, 2009.

\bibitem[BL94]{BarrettLiu:1994}
John~W. Barrett and W.~B. Liu.
\newblock Finite element approximation of the parabolic {$p$}-{L}aplacian.
\newblock {\em SIAM J. Numer. Anal.}, 31(2):413--428, 1994.

\bibitem[Bra01]{Braess:2001}
Dietrich Braess.
\newblock {\em Finite elements}.
\newblock Cambridge University Press, Cambridge, second edition, 2001.
\newblock Theory, fast solvers, and applications in solid mechanics, Translated
  from the 1992 German edition by Larry L. Schumaker.

\bibitem[BS94]{Brenner:1994}
Susanne~C. Brenner and L.~Ridgway Scott.
\newblock {\em The mathematical theory of finite element methods}.
\newblock Springer-Verlag, New York, 1994.

\bibitem[Cia78]{Ciarlet:1978}
Philippe~G. Ciarlet.
\newblock {\em The finite element method for elliptic problems}.
\newblock North-Holland Publishing Co., Amsterdam, 1978.
\newblock Studies in Mathematics and its Applications, Vol. 4.

\bibitem[CMKO11]{ChristiansenMunthe-KaasOwren:2011}
Snorre~H. Christiansen, Hans~Z. Munthe-Kaas, and Brynjulf Owren.
\newblock Topics in structure-preserving discretization.
\newblock {\em Acta Numer.}, 20:1--119, 2011.

\bibitem[CMR09]{CangianiManziniRusso:2009}
Andrea Cangiani, Gianmarco Manzini, and Alessandro Russo.
\newblock Convergence analysis of the mimetic finite difference method for
  elliptic problems.
\newblock {\em SIAM J. Numer. Anal.}, 47(4):2612--2637, 2009.

\bibitem[Dor01]{Dorodnitsyn:2001}
Vladimir Dorodnitsyn.
\newblock Noether-type theorems for difference equations.
\newblock {\em Appl. Numer. Math.}, 39(3-4):307--321, 2001.
\newblock Special issue: Themes in geometric integration.

\bibitem[EG04]{ErnGuermond:2004}
Alexandre Ern and Jean-Luc Guermond.
\newblock {\em Theory and practice of finite elements}, volume 159 of {\em
  Applied Mathematical Sciences}.
\newblock Springer-Verlag, New York, 2004.

\bibitem[Eva98]{Evans:1998}
Lawrence~C. Evans.
\newblock {\em Partial differential equations}, volume~19 of {\em Graduate
  Studies in Mathematics}.
\newblock American Mathematical Society, Providence, RI, 1998.

\bibitem[GH96]{GiaquintaHildebrandt1:1996}
Mariano Giaquinta and Stefan Hildebrandt.
\newblock {\em Calculus of variations. {I}}, volume 310 of {\em Grundlehren der
  Mathematischen Wissenschaften [Fundamental Principles of Mathematical
  Sciences]}.
\newblock Springer-Verlag, Berlin, 1996.
\newblock The Lagrangian formalism.

\bibitem[HLW06]{HairerLubichWanner:2006}
Ernst Hairer, Christian Lubich, and Gerhard Wanner.
\newblock {\em Geometric numerical integration}, volume~31 of {\em Springer
  Series in Computational Mathematics}.
\newblock Springer-Verlag, Berlin, second edition, 2006.
\newblock Structure-preserving algorithms for ordinary differential equations.

\bibitem[HM04]{HydonMansfield:2004}
Peter~E. Hydon and Elizabeth~L. Mansfield.
\newblock A variational complex for difference equations.
\newblock {\em Found. Comput. Math.}, 4(2):187--217, 2004.

\bibitem[Lie71]{Lie:1971}
Sophus Lie.
\newblock {\em Vorlesungen \"uber continuierliche {G}ruppen mit {G}eometrischen
  und anderen {A}nwendungen}.
\newblock Chelsea Publishing Co., Bronx, N.Y., 1971.
\newblock Bearbeitet und herausgegeben von Georg Scheffers, Nachdruck der
  Auflage des Jahres 1893.

\bibitem[Man06]{Mansfield:2006}
Elizabeth~L. Mansfield.
\newblock Noether's theorem for smooth, difference and finite element systems.
\newblock In {\em Foundations of computational mathematics, {S}antander 2005},
  volume 331 of {\em London Math. Soc. Lecture Note Ser.}, pages 230--254.
  Cambridge Univ. Press, Cambridge, 2006.

\bibitem[MP12]{MansfieldPryer:2012}
Elizabeth~L. Mansfield and Tristan. Pryer.
\newblock A noether type theorem for discontinuous galerkin finite element
  schemes.
\newblock {\em In preparation.}, 2012.

\bibitem[Noe71]{Noether:1971}
Emmy Noether.
\newblock Invariant variation problems.
\newblock {\em Transport Theory Statist. Phys.}, 1(3):186--207, 1971.
\newblock Translated from the German (Nachr. Akad. Wiss. G{\"o}ttingen
  Math.-Phys. Kl. II 1918, 235--257).

\bibitem[Olv93]{Olver:1993}
Peter~J. Olver.
\newblock {\em Applications of {L}ie groups to differential equations}, volume
  107 of {\em Graduate Texts in Mathematics}.
\newblock Springer-Verlag, New York, second edition, 1993.

\bibitem[RT77]{RaviartThomas:1977}
P.-A. Raviart and J.~M. Thomas.
\newblock A mixed finite element method for 2nd order elliptic problems.
\newblock In {\em Mathematical aspects of finite element methods ({P}roc.
  {C}onf., {C}onsiglio {N}az. delle {R}icerche ({C}.{N}.{R}.), {R}ome, 1975)},
  pages 292--315. Lecture Notes in Math., Vol. 606. Springer, Berlin, 1977.

\bibitem[SS05]{alberta}
Alfred Schmidt and Kunibert~G. Siebert.
\newblock {\em Design of adaptive finite element software}, volume~42 of {\em
  Lecture Notes in Computational Science and Engineering}.
\newblock Springer-Verlag, Berlin, 2005.
\newblock The finite element toolbox ALBERTA, With 1 CD-ROM (Unix/Linux).

\end{thebibliography}

\end{document}